# Covariance in Non-Commutative Algebra

Aleks Kleyn

ABSTRACT. Consider vector space over non-commutative division algebra. Set of automorphisms of this vector space is group $GL$. Group $GL$ acts on the set of bases of vector space (basis manifold) single transitive and generates active representation. Twin representation on basis manifold is called passive representation. There is no automorphism associated with passive transformation. However passive transformation generates transformation of coordinates of vector with respect to basis. If we consider homomorphism of vector space $V$ into vector space $W$, then we can learn how passive transformation in vector space $V$ generates transformation of coordinates of vector in vector space $W$. Vector in vector space $W$ is called geometric object in vector space $V$. Covariance principle states that geometric object does not depend on the choice of basis. I considered transformation of coordinates of vector and polylinear map.

## Contents

## 1. Preface

Consider rubber tube directed vertically. The air that moves inside the pipe changes the surface of the tube such way that the horizontal section has the form

$$\frac{x^2}{a^2} + \frac{y^2}{b^2} = 1 \qquad (1.1)$$

Aleks_Kleyn@MailAPS.org.
  http://AleksKleyn.dyndns-home.com:4080/   http://arxiv.org/a/kleyn_a_1.
  http://AleksKleyn.blogspot.com/.

where parameters $a$ and $b$ are functions of coordinate $z$ and time $t$ and are solution of the system of differential equations

(1.2)
$$\frac{\partial a}{\partial z} = a_1 x + b_1 y \quad \frac{\partial b}{\partial z} = c_1 x + d_1 y$$
$$\frac{\partial a}{\partial t} = a_2 x + b_2 y \quad \frac{\partial b}{\partial t} = c_2 x + d_2 y$$

However, we can measure distance in meters or centimeters; as well we can measure time in minutes or seconds. We can change direction of axes $x$, $y$ in a plane perpendicular to the axis $z$. These coordinate transformations will be reflected in the notation of the equation and its solution. However, the solution to a problem does not depend on the form of this solution in coordinate representation.

Consider more simple example. Curve of the second order on the plane has a very complex representation. However, we can always choose coordinates with respect to which the curve takes the canonical form and is an ellipse, hyperbole, parabola or a pair of straight lines.

Symmetry consisting in the independence of the solution from the choice of coordinates is called covariance.

Tensor algebra is the expression of covariance in vector space over field. Our goal is to consider the theory of covariance in vector space over $D$-algebra. Although this theory is more complex than theory of covariance in vector space over field, both theories share common features.

## 2. Product of matrices in non-commutative algebra

Before we begin to consider vector space over non-commutative division algebra $A$, I recall that there are two operations of product of matrices with entries from algebra $A$.

**Definition 2.1.** *Let the nubmer of columns of the matrix $a$ equal the number of rows of the matrix $b$. $_*^*$-**product** of matrices $a$ and $b$ has form*

(2.1)
$$a_*{}^*b = \left(a_k^i b_j^k\right)$$

(2.2)
$$(a_*{}^*b)_j^i = a_k^i b_j^k$$

(2.3)
$$\begin{pmatrix} a_1^1 & \ldots & a_p^1 \\ \ldots & \ldots & \ldots \\ a_1^n & \ldots & a_p^n \end{pmatrix} {}_*^* \begin{pmatrix} b_1^1 & \ldots & b_m^1 \\ \ldots & \ldots & \ldots \\ b_1^p & \ldots & b_m^p \end{pmatrix} = \begin{pmatrix} a_k^1 b_1^k & \ldots & a_k^1 b_m^k \\ \ldots & \ldots & \ldots \\ a_k^n b_1^k & \ldots & a_k^n b_m^k \end{pmatrix}$$
$$= \begin{pmatrix} (a_*{}^*b)_1^1 & \ldots & (a_*{}^*b)_m^1 \\ \ldots & \ldots & \ldots \\ (a_*{}^*b)_1^n & \ldots & (a_*{}^*b)_m^n \end{pmatrix}$$

*$_*^*$-product can be expressed as product of a row of the matrix $a$ over a column of the matrix $b$.* □

**Definition 2.2.** *Let the nubmer of rows of the matrix $a$ equal the number of columns of the matrix $b$. $^*{}_*$-**product** of matrices $a$ and $b$ has form*

$$(2.4) \qquad a^*{}_*b = \left(a_i^k b_k^j\right)$$

$$(2.5) \qquad (a^*{}_*b)_j^i = a_i^k b_k^j$$

$$(2.6) \quad \begin{pmatrix} a_1^1 & \dots & a_m^1 \\ \dots & \dots & \dots \\ a_1^p & \dots & a_m^p \end{pmatrix} {}^*{}_* \begin{pmatrix} b_1^1 & \dots & b_p^1 \\ \dots & \dots & \dots \\ b_1^n & \dots & b_p^n \end{pmatrix} = \begin{pmatrix} a_1^k b_k^1 & \dots & a_m^k b_k^1 \\ \dots & \dots & \dots \\ a_1^k b_k^n & \dots & a_m^k b_k^n \end{pmatrix}$$
$$= \begin{pmatrix} (a^*{}_*b)_1^1 & \dots & (a^*{}_*b)_m^1 \\ \dots & \dots & \dots \\ (a^*{}_*b)_1^n & \dots & (a^*{}_*b)_m^n \end{pmatrix}$$

$^*{}_*$-*product can be expressed as product of a column of the matrix $a$ over a row of the matrix $b$.* □

## 3. Homomorphism of Left $A$-Vector Space

Let $A$ be non-commutative division $D$-algebra. Then we can consider left or right $A$-vector space. We may confine ourselves to considering of left $A$-vector space, because statements for right $A$-vector space are similar.

We will start by considering of the set of homomorphisms of left $A$-vector space. Homomorphism is a map that preserves the structure of algebra. In particular, the theorem 3.2 is true.

**Definition 3.1.** *Let diagram of representations*

$$(3.1) \quad \begin{array}{c} A \xrightarrow{g_{1.23}} A \xrightarrow{g_{1.34}} V_1 \\ {}_{g_{1.12}}\nwarrow \quad \uparrow {}_{g_{1.12}} \quad \nearrow {}_{g_{1.14}} \\ D \end{array} \qquad \begin{array}{l} g_{12}(d): a \to d\,a \\ g_{23}(v): w \to C(w,v) \\ C \in \mathcal{L}(A^2 \to A) \\ g_{3,4}(a): v \to a\,v \\ g_{1,4}(d): v \to d\,v \end{array}$$

*describe left $A$-module $V_1$. Let diagram of representations*

$$(3.2) \quad \begin{array}{c} A \xrightarrow{g_{2.23}} A \xrightarrow{g_{2.34}} V_2 \\ {}_{g_{2.12}}\nwarrow \quad \uparrow {}_{g_{2.12}} \quad \nearrow {}_{g_{2.14}} \\ D \end{array} \qquad \begin{array}{l} g_{12}(d): a \to d\,a \\ g_{23}(v): w \to C(w,v) \\ C \in \mathcal{L}(A^2 \to A) \\ g_{3,4}(a): v \to a\,v \\ g_{1,4}(d): v \to d\,v \end{array}$$

*describe left $A$-module $V_2$. Morphism*

$$(3.3) \qquad f: V_1 \to V_2$$

of diagram of representations (3.1) into diagram of representations (3.2) *is called* **homomorphism** *of left $A$-module $V_1$ into left $A$-module $V_2$. Let us denote* $\text{Hom}(D; A; *V_1 \to *V_2)$ *set of homomorphisms of left $A$-module $V_1$ into left $A$-module $V_2$.* □

We will use notation
$$f \circ a = f(a)$$
for image of homomorphism $f$.

**Theorem 3.2.** *The homomorphism*

(3.3)  $f : V_1 \to V_2$

*of left $A$-module $V_1$ into left $A$-module $V_2$ satisfies following equalities*

(3.4) $$f \circ (u + v) = f \circ u + f \circ v$$

(3.5) $$f \circ (av) = a(f \circ v)$$

$$a \in A \quad u, v \in V_1$$

PROOF. The theorem follows from the theorem [3]-10.3.2. □

**Definition 3.3.** *Homomorphism*
$$f : V \to W$$
*is called*[3.1]**isomorphism** *between left $A$-vector spaces $V$ and $W$, if correspondence $f^{-1}$ is homomorphism. A homomorphism*
$$f : V \to V$$
*in which source and target are the same left $A$-vector space is called* **endomorphism**. *Endomorphism*
$$f : V \to V$$
*of left $A$-vector space $V$ is called* **automorphism**, *if correspondence $f^{-1}$ is endomorphism.* □

**Theorem 3.4.** *The set $GL(V)$ of automorphisms of left $A$-vector space $V$ is group.*

Although the theorem 3.2 is true and can serve as a definition of homomorphism, this theorem does not give the full picture.

Suppose we are studying a physical phenomenon. We can give a qualitative picture of this phenomenon. However, if we want to perform experiment, we have to choose a set of measuring tools and to calculate expected measurements. The set of measuring tools is a basis in the space of measurements; and expected measurements are coordinates of experiment in the space of measurements.

We will extend this remark to the theorem 3.2. To understand the structure of homomorphism and answer the question how big is the set of homomorphisms, we have to choose a basis of $A$-vector space and consider the change of coordinates of vector with respect to selected basis.

---

[3.1] I follow the definition on page [4]-49.

There are various conventions for numbering vectors of the basis and corresponding conventions for the numbering of coordinates of a vector. For instance, in the theorem 3.5, we consider coordinates of vector in left $A$-vector space of columns.

**Theorem 3.5.** *If we write vectors of basis $\overline{\overline{e}}$ as row of matrix*

$$(3.6) \qquad e = \begin{pmatrix} e_1 & ... & e_n \end{pmatrix}$$

*and coordinates of vector $\overline{w} = w^i e_i$ with respect to basis $\overline{\overline{e}}$ as column of matrix*

$$(3.7) \qquad w = \begin{pmatrix} w^1 \\ ... \\ w^n \end{pmatrix}$$

*then we can represent the vector $\overline{w}$ as $^*_*$-product of matrices*

$$(3.8) \qquad \overline{w} = w {^*_*} e = \begin{pmatrix} w^1 \\ ... \\ w^n \end{pmatrix} {^*_*} \begin{pmatrix} e_1 & ... & e_n \end{pmatrix} = w^i e_i$$

If the basis $\overline{\overline{e}}$ is given, then

$$(3.9) \qquad V_* = \{ v : \overline{v} = v {^*_*} e, \overline{v} \in V \}$$

is the set of coordinates of vectors $\overline{v} \in V$. The set $V_*$ is left $A$-vector space and is isomorphic to left $A$-vector space $V$. Therefore, the set $V_*$ does not depend on choice of basis $\overline{\overline{e}}$.

If left $A$-vector space $V$ has dimension $n$, then left $A$-vector space $V_*$ is isomorphic to direct sum of $n$ copies of $D$-algebra $A$. So in this case we put

$$(3.10) \qquad V_* = n {^*_*} A$$

In the theorem 3.6, we consider coordinates of vector in left $A$-vector space of rows.

**Theorem 3.6.** *If we write vectors of basis $\overline{\overline{e}}$ as column of matrix*

$$(3.11) \qquad e = \begin{pmatrix} e^1 \\ ... \\ e^n \end{pmatrix}$$

*and coordinates of vector $\overline{w} = w_i e^i$ with respect to basis $\overline{\overline{e}}$ as row of matrix*

$$(3.12) \qquad w = \begin{pmatrix} w_1 & ... & w_n \end{pmatrix}$$

*then we can represent the vector $\overline{w}$ as $_*^*$-product of matrices*

$$(3.13) \qquad \overline{w} = w {_*^*} e = \begin{pmatrix} w_1 & ... & w_n \end{pmatrix} {_*^*} \begin{pmatrix} e^1 \\ ... \\ e^n \end{pmatrix} = w_i e^i$$

We considered main formats of presentation of coordinates of vector of left $A$-vector space. It is evident that other formats of presentation of coordinates of vector are also possible. For instance, we may consider the set of $n \times m$ matrices as left $A$-vector space. The same way we may consider main formats of presentation of coordinates of vector of right $A$-vector space.

Statements for various forms of representation of coordinates of vector are similar. Therefore, we will focus our attention on left $A$-vector space of columns.

**Theorem 3.7.** *The homomorphism* [3.2]

$$(3.14) \qquad \overline{f}: V_1 \to V_2$$

*of left $A$-vector space $V_1$ into left $A$-vector space $V_2$ has presentation*

$$(3.15) \qquad w = v^*{}_* f$$

$$(3.16) \qquad \overline{f} \circ (v^i e_{V_1 i}) = v^i f_i^k e_{V_2 k}$$

$$(3.17) \qquad \overline{f} \circ (v^*{}_* e_1) = v^*{}_* f^*{}_* e_2$$

*relative to selected bases. Here*

- $v$ *is coordinate matrix of $V_1$-number $\overline{v}$ relative the basis $\overline{\overline{e}}_{V_1}$*

$$(3.18) \qquad \overline{v} = v^*{}_* e_{V_1}$$

- $w$ *is coordinate matrix of $V_2$-number*

$$(3.19) \qquad \overline{w} = \overline{f} \circ \overline{v}$$

*relative the basis $\overline{\overline{e}}_{V_2}$*

$$(3.20) \qquad \overline{w} = w^*{}_* e_{V_2}$$

- $f$ *is coordinate matrix of set of $V_2$-numbers $(\overline{f} \circ e_{V_1 i}, i \in I)$ relative the basis $\overline{\overline{e}}_{V_2}$.*

*The matrix $f$ is unique and is called* **matrix of homomorphism $\overline{f}$** *relative bases $\overline{\overline{e}}_1$, $\overline{\overline{e}}_2$.*

PROOF. The theorem follows from the theorem [3]-10.3.3. □

The converse theorem is also true.

**Theorem 3.9.** *Let*

$$f = (f_j^i, i \in I, j \in J)$$

**Theorem 3.8.** *The homomorphism* [3.3]

$$(3.21) \qquad \overline{f}: V_1 \to V_2$$

*of left $A$-vector space $V_1$ into left $A$-vector space $V_2$ has presentation*

$$(3.22) \qquad w = v_*{}^* f$$

$$(3.23) \qquad \overline{f} \circ (v_i e_{V_1}^i) = v_i f_k^i e_{V_2}^k$$

$$(3.24) \qquad \overline{f} \circ (v_*{}^* e_1) = v_*{}^* f_*{}^* e_2$$

*relative to selected bases. Here*

- $v$ *is coordinate matrix of $V_1$-number $\overline{v}$ relative the basis $\overline{\overline{e}}_{V_1}$*

$$(3.25) \qquad \overline{v} = v_*{}^* e_{V_1}$$

- $w$ *is coordinate matrix of $V_2$-number*

$$(3.26) \qquad \overline{w} = \overline{f} \circ \overline{v}$$

*relative the basis $\overline{\overline{e}}_{V_2}$*

$$(3.27) \qquad \overline{w} = w_*{}^* e_{V_2}$$

- $f$ *is coordinate matrix of set of $V_2$-numbers $(\overline{f} \circ e_{V_1}^i, i \in I)$ relative the basis $\overline{\overline{e}}_{V_2}$.*

*The matrix $f$ is unique and is called* **matrix of homomorphism $\overline{f}$** *relative bases $\overline{\overline{e}}_1$, $\overline{\overline{e}}_2$.*

PROOF. The theorem follows from the theorem [3]-10.3.4. □

**Theorem 3.10.** *Let*

$$f = (f_i^j, i \in I, j \in J)$$

---

[3.2] In theorems 3.7, 3.9, we use the following convention. Let the set of vectors $\overline{\overline{e}}_1 = (e_{1i}, i \in I)$ be a basis of left $A$-vector space $V_1$. Let the set of vectors $\overline{\overline{e}}_2 = (e_{2j}, j \in J)$ be a basis of left $A$-vector space $V_2$.

[3.3] In theorems 3.8, 3.10, we use the following convention. Let the set of vectors $\overline{\overline{e}}_1 = (e_1^i, i \in I)$ be a basis of left $A$-vector space $V_1$. Let the set of vectors $\overline{\overline{e}}_2 = (e_2^j, j \in J)$ be a basis of left $A$-vector space $V_2$.

| | |
|---|---|
| be matrix of $A$-numbers. The map[3.2] $$(3.14) \quad \overline{f} : V_1 \to V_2$$ defined by the equality $$(3.17) \quad \overline{f} \circ (v^*{}_* e_1) = v^*{}_* f^*{}_* e_2$$ is homomorphism of left $A$-vector space of columns. The homomorphism (3.14) which has the given matrix $f$ is unique. | be matrix of $A$-numbers. The map[3.3] $$(3.21) \quad \overline{f} : V_1 \to V_2$$ defined by the equality $$(3.24) \quad \overline{f} \circ (v_*{}^* e_1) = v_*{}^* f_*{}^* e_2$$ is homomorphism of left $A$-vector space of rows. The homomorphism (3.21) which has the given matrix $f$ is unique. |

PROOF. The theorem follows from the theorem [3]-10.3.5. □

PROOF. The theorem follows from the theorem [3]-10.3.6. □

Therefore, if we choose the basis $\overline{\overline{e}} = (e_1, ..., e_n)$ of left $A$-vector space $V$ of columns, then we can identify the set of endomorphisms of left $A$-vector space $V$ and the set of $n \times n$ matrices. We will see this relation between endomorphisms and matrices deeper when we consider automorphisms.

| | |
|---|---|
| **Theorem 3.11.** *Let $V$ be a left $A$-vector space of columns and $\overline{\overline{e}}$ be basis of left $A$-vector space $V$. Any automorphism $\overline{f}$ of left $A$-vector space $V$ has form* $$(3.28) \qquad v' = v^*{}_* f$$ *where $f$ is a $^*{}_*$-nonsingular matrix. Matrices of automorphisms of left $A$-vector space $V$ of columns form a group $GL(V_*)$ isomorphic to group $GL(V)$. Automorphisms of left $A$-vector space of columns form a right-side linear effective representation* $$(3.29) \qquad GL(V_*) \relbar\joinrel\ast\joinrel\rightarrow V_*$$ *of the group $GL(V_*)$ in left $A$-vector space $V_*$.* | **Theorem 3.12.** *Let $V$ be a left $A$-vector space of rows and $\overline{\overline{e}}$ be basis of left $A$-vector space $V$. Any automorphism $\overline{f}$ of left $A$-vector space $V$ has form* $$(3.30) \qquad v' = v_*{}^* f$$ *where $f$ is a $_*{}^*$-nonsingular matrix. Matrices of automorphisms of left $A$-vector space $V$ of rows form a group $GL(V_*)$ isomorphic to group $GL(V)$. Automorphisms of left $A$-vector space of rows form a right-side linear effective representation* $$(3.31) \qquad GL(V_*) \relbar\joinrel\ast\joinrel\rightarrow V_*$$ *of the group $GL(V_*)$ in left $A$-vector space $V_*$.* |

PROOF. The theorem follows from the theorem [3]-12.1.9, [3]-12.1.11, [3]-12.2.1. □

PROOF. The theorem follows from the theorem [3]-12.1.10, [3]-12.1.12, [3]-12.2.2. □

Theorems 3.14, 3.16 show relationship between $_*{}^*$-nonsingular matrices, bases and automorphisms of left $A$-vector space of rows.

| | |
|---|---|
| **Theorem 3.13.** *Let $V$ be a left $A$-vector space of columns. The coordinate matrix of basis $\overline{\overline{g}}$ relative basis $\overline{\overline{e}}$ of left $A$-vector space $V$ is $^*{}_*$-nonsingular matrix.* | **Theorem 3.14.** *Let $V$ be a left $A$-vector space of rows. The coordinate matrix of basis $\overline{\overline{g}}$ relative basis $\overline{\overline{e}}$ of left $A$-vector space $V$ is $_*{}^*$-nonsingular matrix.* |

PROOF. The theorem follows from the theorem [3]-12.1.7. □

PROOF. The theorem follows from the theorem [3]-12.1.8. □

| **Theorem 3.15.** *Automorphism a acting on each vector of basis of left A-vector space of columns maps a basis into another basis.* | **Theorem 3.16.** *Automorphism a acting on each vector of basis of left A-vector space of rows maps a basis into another basis.* |

PROOF. The theorem follows from the theorem [3]-12.2.3. □

PROOF. The theorem follows from the theorem [3]-12.2.4. □

Thus we can extend a right-side linear $GL(V_*)$-representation in left $A$-vector space $V_*$ to the set of bases of left $A$-vector space $V$. Transformation of this right-side representation on the set of bases of left $A$-vector space $V$ is called **active transformation** because the homomorphism of the left $A$-vector space induced this transformation (See also definition in the section [2]-14.1-3 as well the definition on the page [1]-214).

According to definition we write the action of the transformation $a \in GL(V_*)$ on the basis $\overline{\overline{e}}$ as $\overline{\overline{e}}^*{}_*a$. Consider the equality

$$(3.32) \qquad v^*{}_*\overline{\overline{e}}^*{}_*a = v^*{}_*\overline{\overline{e}}^*{}_*a$$

The expression $\overline{\overline{e}}^*{}_*a$ on the left side of the equality (3.32) is image of basis $\overline{\overline{e}}$ with respect to active transformation $a$. The expression $v^*{}_*\overline{\overline{e}}$ on the right side of the equality (3.32) is expansion of vector $\overline{v}$ with respect to basis $\overline{\overline{e}}$. Therefore, the expression on the right side of the equality (3.32) is image of vector $\overline{v}$ with respect to endomorphism $a$ and the expression on the left side of the equality

$(3.32) \quad v^*{}_*\overline{\overline{e}}^*{}_*a = v^*{}_*\overline{\overline{e}}^*{}_*a$

is expansion of image of vector $\overline{v}$ with respect to image of basis $\overline{\overline{e}}$. Therefore, from the equality (3.32) it follows that endomorphism $a$ of left $A$-vector space and corresponding active transformation $a$ act synchronously and coordinates $a \circ \overline{v}$ of image of the vector $\overline{v}$ with respect to the image $\overline{\overline{e}}^*{}_*a$ of the basis $\overline{\overline{e}}$ are the same as coordinates of the vector $\overline{v}$ with respect to the basis $\overline{\overline{e}}$.

According to definition we write the action of the transformation $a \in GL(V_*)$ on the basis $\overline{\overline{e}}$ as $\overline{\overline{e}}_*{}^*a$. Consider the equality

$$(3.33) \qquad v_*{}^*\overline{\overline{e}}_*{}^*a = v_*{}^*\overline{\overline{e}}_*{}^*a$$

The expression $\overline{\overline{e}}_*{}^*a$ on the left side of the equality (3.33) is image of basis $\overline{\overline{e}}$ with respect to active transformation $a$. The expression $v_*{}^*\overline{\overline{e}}$ on the right side of the equality (3.33) is expansion of vector $\overline{v}$ with respect to basis $\overline{\overline{e}}$. Therefore, the expression on the right side of the equality (3.33) is image of vector $\overline{v}$ with respect to endomorphism $a$ and the expression on the left side of the equality

$(3.33) \quad v_*{}^*\overline{\overline{e}}_*{}^*a = v_*{}^*\overline{\overline{e}}_*{}^*a$

is expansion of image of vector $\overline{v}$ with respect to image of basis $\overline{\overline{e}}$. Therefore, from the equality (3.33) it follows that endomorphism $a$ of left $A$-vector space and corresponding active transformation $a$ act synchronously and coordinates $a \circ \overline{v}$ of image of the vector $\overline{v}$ with respect to the image $\overline{\overline{e}}_*{}^*a$ of the basis $\overline{\overline{e}}$ are the same as coordinates of the vector $\overline{v}$ with respect to the basis $\overline{\overline{e}}$.

## 4. Passive Transformation

Let us define an additional structure on left $A$-vector space $V$. Then not every automorphism keeps properties of the selected structure. For imstance, if we introduce norm in left $A$-vector space $V$, then we are interested in automorphisms which preserve the norm of the vector.

**Definition 4.1.** *Normal subgroup $G(V)$ of the group $GL(V)$ such that subgroup $G(V)$ generates automorphisms which hold properties of the selected structure is called* **symmetry group**.

*Without loss of generality we identify element $g$ of group $G(V)$ with corresponding transformation of representation and write its action on vector $v \in V$ as $v^*{}_*g$.* □

**Definition 4.2.** *Normal subgroup $G(V)$ of the group $GL(V)$ such that subgroup $G(V)$ generates automorphisms which hold properties of the selected structure is called* **symmetry group**.

*Without loss of generality we identify element $g$ of group $G(V)$ with corresponding transformation of representation and write its action on vector $v \in V$ as $v_*{}^*g$.* □

**Definition 4.3.** *The right-side representation of group $G(V)$ in the set of bases of left $A$-vector space $V$ is called* **active left-side $G$-representation**. □

If the basis $\bar{\bar{e}}$ is given, then we can identify the automorphism of left $A$-vector space $V$ and its coordinates with respect to the basis $\bar{\bar{e}}$. The set $G(V_*)$ of coordinates of automorphisms with respect to the basis $\bar{\bar{e}}$ is group isomorphic to the group $G(V)$.

Not every two bases can be mapped by a transformation from the symmetry group because not every nonsingular linear transformation belongs to the representation of group $G(V)$. Therefore, we can represent the set of bases as a union of orbits of group $G(V)$. In particular, if the basis $\bar{\bar{e}} \in G(V)$, then the group $G(V)$ is orbit of the basis $\bar{\bar{e}}$.

**Definition 4.4.** *We call orbit $\bar{\bar{e}}^*{}_*G(V)$ of the selected basis $\bar{\bar{e}}$ the* **basis manifold** *of left $A$-vector space $V$ of columns.* □

**Definition 4.5.** *We call orbit $\bar{\bar{e}}_*{}^*G(V)$ of the selected basis $\bar{\bar{e}}$ the* **basis manifold** *of left $A$-vector space $V$ of rows.* □

**Theorem 4.6.** *Active right-side $G(V)$-representation on basis manifold is single transitive representation.*

PROOF. The theorem follows from the theorem [3]-12.2.12. □

**Theorem 4.7.** *On the basis manifold $\bar{\bar{e}}^*{}_*G(V)$ of left $A$-vector space of columns, there exists single transitive left-side $G(V)$-representation, commuting with active.*

**Theorem 4.8.** *On the basis manifold $\bar{\bar{e}}_*{}^*G(V)$ of left $A$-vector space of rows, there exists single transitive left-side $G(V)$-representation, commuting with active.*

PROOF. The theorem follows from the theorem [3]-12.2.13. □

PROOF. The theorem follows from the theorem [3]-12.2.14. □

Transformation of left-side $G(V)$-representation is different from an active transformation and cannot be reduced to transformation of space $V$.

**Definition 4.9.** *A transformation of left-side $G(V)$-representation is called* **passive transformation** *of basis manifold $\overline{\overline{e}}{}^*{}_*G(V)$ of left $A$-vector space of columns, and the left-side $G(V)$-representation is called* **passive left-side $G(V)$-representation**. *According to the definition we write the passive transformation of basis $\overline{\overline{e}}$ defined by element $a \in G(V)$ as $a^*{}_*\overline{\overline{e}}$.* □

**Definition 4.10.** *A transformation of left-side $G(V)$-representation is called* **passive transformation** *of basis manifold $\overline{\overline{e}}_*{}^*G(V)$ of left $A$-vector space of rows, and the left-side $G(V)$-representation is called* **passive left-side $G(V)$-representation**. *According to the definition we write the passive transformation of basis $\overline{\overline{e}}$ defined by element $a \in G(V)$ as $a_*{}^*\overline{\overline{e}}$.* □

## 5. Geometric Object

An active transformation changes bases and vectors uniformly and coordinates of vector relative basis do not change. A passive transformation changes only the basis and it leads to transformation of coordinates of vector relative to basis.

We consider transformation of coordinates of vector in the theorem 5.1.

We consider transformation of coordinates of vector in the theorem 5.2.

**Theorem 5.1.** *Let $V$ be a left $A$-vector space of columns. Let passive transformation $g \in G(V_*)$ map the basis $\overline{\overline{e}}_1$ into the basis $\overline{\overline{e}}_2$*

$$(5.1) \qquad \overline{\overline{e}}_2 = g^*{}_*\overline{\overline{e}}_1$$

*Let*

$$(5.2) \qquad v_i = \begin{pmatrix} v_i^{\textbf{\textit{1}}} \\ ... \\ v_i^{\textbf{\textit{n}}} \end{pmatrix}$$

*be matrix of coordinates of the vector $\overline{v}$ with respect to the basis $\overline{\overline{e}}_i$, $i = 1, 2$. Coordinate transformations*

$$(5.3) \qquad v_1 = v_2{}^*{}_*g$$

$$(5.4) \qquad v_2 = v_1{}^*{}_*g^{-1^*{}_*}$$

*do not depend on vector $\overline{v}$ or basis $\overline{\overline{e}}$, but is defined only by coordinates of vector $\overline{v}$ relative to basis $\overline{\overline{e}}$.*

Proof. The theorem follows from the theorem [3]-12.3.1. □

Let $V$, $W$ be left $A$-vector spaces of columns and $G(V_*)$ be symmetry group of left $A$-vector space $V$. Homomor-

**Theorem 5.2.** *Let $V$ be a left $A$-vector space of rows. Let passive transformation $g \in G(V_*)$ map the basis $\overline{\overline{e}}_1$ into the basis $\overline{\overline{e}}_2$*

$$(5.5) \qquad \overline{\overline{e}}_2 = g_*{}^*\overline{\overline{e}}_1$$

*Let*

$$v_i = \begin{pmatrix} v_{i\textbf{\textit{1}}} & ... & v_{i\textbf{\textit{n}}} \end{pmatrix}$$

*be matrix of coordinates of the vector $\overline{v}$ with respect to the basis $\overline{\overline{e}}_i$, $i = 1, 2$. Coordinate transformations*

$$(5.6) \qquad v_1 = v_{2*}{}^*g$$

$$(5.7) \qquad v_2 = v_{1*}{}^*g^{-1_*{}^*}$$

*do not depend on vector $\overline{v}$ or basis $\overline{\overline{e}}$, but is defined only by coordinates of vector $\overline{v}$ relative to basis $\overline{\overline{e}}$.*

Proof. The theorem follows from the theorem [3]-12.3.2. □

Let $V$, $W$ be left $A$-vector spaces of rows and $G(V_*)$ be symmetry group of left $A$-vector space $V$. Homomorphism

$$(5.13) \qquad F : G(V_*) \to GL(W_*)$$

phism

$$(5.8) \qquad F : G(V_*) \to GL(W_*)$$

maps passive transformation $g \in G(V_*)$

$$(5.9) \qquad \overline{\overline{e}}_{V2} = g^*{}_*\overline{\overline{e}}_{V1}$$

of left $A$-vector spaces $V$ into passive transformation $F(g) \in GL(W_*)$

$$(5.10) \qquad \overline{\overline{e}}_{W2} = F(g)^*{}_*\overline{\overline{e}}_{W1}$$

of left $A$-vector spaces $W$.

Then coordinate transformation in left $A$-vector space $W$ gets form

$$(5.11) \qquad w_2 = w_1{}^*{}_* F(g)^{-1^*{}_*}$$

Therefore, the map

$$(5.12) \qquad F_1(g) = F(g)^{-1^*{}_*}$$

is right-side representation

$$F_1 : G(V_*) \longrightarrow{*} W_*$$

of group $G(V_*)$ in the set $W_*$.

maps passive transformation $g \in G(V_*)$

$$(5.14) \qquad \overline{\overline{e}}_{V2} = g_*{}^*\overline{\overline{e}}_{V1}$$

of left $A$-vector spaces $V$ into passive transformation $F(g) \in GL(W_*)$

$$(5.15) \qquad \overline{\overline{e}}_{W2} = F(g)_*{}^*\overline{\overline{e}}_{W1}$$

of left $A$-vector spaces $W$.

Then coordinate transformation in left $A$-vector space $W$ gets form

$$(5.16) \qquad w_2 = w_1{}_*{}^* F(g)^{-1_*{}^*}$$

Therefore, the map

$$(5.17) \qquad F_1(g) = F(g)^{-1_*{}^*}$$

is right-side representation

$$F_1 : G(V_*) \longrightarrow{*} W_*$$

of group $G(V_*)$ in the set $W_*$.

**Definition 5.3.** *Orbit*

$$(5.18) \qquad \begin{aligned} & \mathcal{O}(V, W, \overline{\overline{e}}_V, w) \\ & = (w^*{}_* F(G)^{-1^*{}_*}, G^*{}_* \overline{\overline{e}}_V) \end{aligned}$$

*of representation $F_1$ is called* **coordinate manifold of geometric object** *in left $A$-vector space $V$ of columns. For any basis*

$(5.9) \quad \overline{\overline{e}}_{V2} = g^*{}_*\overline{\overline{e}}_{V1}$

*corresponding point*

$(5.11) \quad w_2 = w_1{}^*{}_* F(g)^{-1^*{}_*}$

*of orbit defines* **coordinates of geometric object** *in coordinate left $A$-vector space relative basis $\overline{\overline{e}}_{V2}$.* □

**Definition 5.4.** *Orbit*

$$(5.19) \qquad \begin{aligned} & \mathcal{O}(V, W, \overline{\overline{e}}_V, w) \\ & = (w_*{}^* F(G)^{-1_*{}^*}, G_*{}^* \overline{\overline{e}}_V) \end{aligned}$$

*of representation $F_1$ is called* **coordinate manifold of geometric object** *in left $A$-vector space $V$ of rows. For any basis*

$(5.14) \quad \overline{\overline{e}}_{V2} = g_*{}^*\overline{\overline{e}}_{V1}$

*corresponding point*

$(5.16) \quad w_2 = w_1{}_*{}^* F(g)^{-1_*{}^*}$

*of orbit defines* **coordinates of geometric object** *in coordinate left $A$-vector space relative basis $\overline{\overline{e}}_{V2}$.* □

**Definition 5.5.** *Let us say the coordinates $w_1$ of vector $\overline{w}$ with respect to the basis $\overline{\overline{e}}_{W1}$ are given. The set of vectors*

$$(5.20) \qquad \overline{w}_2 = w_2{}^*{}_* e_{W2}$$

*is called* **geometric object** *defined in left $A$-vector space $V$ of columns. For any basis $\overline{\overline{e}}_{W2}$, corresponding point*

$(5.11) \quad w_2 = w_1{}^*{}_* F(g)^{-1^*{}_*}$

**Definition 5.6.** *Let us say the coordinates $w_1$ of vector $\overline{w}$ with respect to the basis $\overline{\overline{e}}_{W1}$ are given. The set of vectors*

$$(5.22) \qquad \overline{w}_2 = w_{2*}{}^* e_{W2}$$

*is called* **geometric object** *defined in left $A$-vector space $V$ of rows. For any basis $\overline{\overline{e}}_{W2}$, corresponding point*

$(5.16) \quad w_2 = w_1{}_*{}^* F(g)^{-1_*{}^*}$

| | |
|---|---|
| of coordinate manifold defines the vector | of coordinate manifold defines the vector |
| (5.21) $\quad \overline{w}_2 = w_2{}^*{}_* e_{W2}$ | (5.23) $\quad \overline{w}_2 = w_{2*}{}^* e_{W2}$ |
| which is called **representative of geometric object** in left $A$-vector space $V$ in basis $\overline{\overline{e}}_{V2}$. $\square$ | which is called **representative of geometric object** in left $A$-vector space $V$ in basis $\overline{\overline{e}}_{V2}$. $\square$ |

**Theorem 5.7. (Principle of covariance).** *Representative of geometric object does not depend on selection of basis $\overline{\overline{e}}_{V2}$.*

## 6. Examples of Geometric Object

According to the theorem 5.1, a vector is geometric object. The set of endomorphisms of left $A$-vector space as well as the set of linear or polylinear maps are $A$-vector space. So we can ask if endomorphism or polylinear map is geometric object.

| | |
|---|---|
| We consider transformation of coordinates of endomorphism of left $A$-vector space of columns in the theorem 6.1. | We consider transformation of coordinates of endomorphism of left $A$-vector space of rows in the theorem 6.2. |
| **Theorem 6.1.** *Let $V$ be a left $A$-vector space of columns and $\overline{\overline{e}}_1$, $\overline{\overline{e}}_2$ be bases in left $A$-vector space $V$. Let $g$ be passive transformation of basis $\overline{\overline{e}}_1$ into basis $\overline{\overline{e}}_2$* $$\overline{\overline{e}}_2 = g^*{}_* \overline{\overline{e}}_1$$ *Let $f$ be endomorphism of left $A$-vector space $V$. Let $f_i$, $i=1,2$, be the matrix of endomorphism $f$ with respect to the basis $\overline{\overline{e}}_i$. Then* (6.1) $\quad f_2 = g^*{}_* f_1{}^*{}_* g^{-1*}{}_*$ | **Theorem 6.2.** *Let $V$ be a left $A$-vector space of rows and $\overline{\overline{e}}_1$, $\overline{\overline{e}}_2$ be bases in left $A$-vector space $V$. Let $g$ be passive transformation of basis $\overline{\overline{e}}_1$ into basis $\overline{\overline{e}}_2$* $$\overline{\overline{e}}_2 = g_*{}^* \overline{\overline{e}}_1$$ *Let $f$ be endomorphism of left $A$-vector space $V$. Let $f_i$, $i=1,2$, be the matrix of endomorphism $f$ with respect to the basis $\overline{\overline{e}}_i$. Then* (6.2) $\quad f_2 = g_*{}^* f_{1*}{}^* g^{-1}{}_*{}^*$ |
| PROOF. The theorem follows from the theorem [3]-12.3.5. $\square$ | PROOF. The theorem follows from the theorem [3]-12.3.6. $\square$ |

According to the theorem 6.1, an endomorphism of vector space is geometric object. Indeed, when we change basis, we see transformation of coordinates of vector and coordinates of endomorphism. However, change of coordinates agreed and image of vector for given endomorphism does not depend on choice of basis.

Let $\overline{\overline{e}}$ be the basis of left vector space $V$ of columns. Linear map

$$a : V \to V$$

has form

$$(6.3) \quad \begin{pmatrix} v^1 \\ \dots \\ v^n \end{pmatrix} = \begin{pmatrix} a_1^1 & \dots & a_n^1 \\ \dots & \dots & \dots \\ a_1^n & \dots & a_n^n \end{pmatrix} \circ^\circ \begin{pmatrix} w^1 \\ \dots \\ w^n \end{pmatrix} = \begin{pmatrix} a_i^1 \circ w^i \\ \dots \\ a_i^n \circ w^i \end{pmatrix}$$

$$a = \begin{pmatrix} a_1^1 & ... & a_n^1 \\ ... & ... & ... \\ a_1^n & ... & a_n^n \end{pmatrix}$$

with respect to the basis $\overline{\overline{e}}$. Since partial map $a_j^i$ is linear map of $D$-algebra $A$, then we can write linear map $a_j^i$ as

(6.4) $\quad a_j^i = a_{j\,s0}^i \otimes a_{j\,s1}^i$

The equality

(6.5) $\quad v^i = (a_{j\,s1}^i \otimes a_{j\,s2}^i) \circ w^j = a_{j\,s1}^i w^j a_{j\,s2}^i$

follows from equalities

(6.3) $\begin{pmatrix} v^1 \\ ... \\ v^n \end{pmatrix} = \begin{pmatrix} a_1^1 & ... & a_n^1 \\ ... & ... & ... \\ a_1^n & ... & a_n^n \end{pmatrix} \circ^{\circ} \begin{pmatrix} w^1 \\ ... \\ w^n \end{pmatrix} = \begin{pmatrix} a_i^1 \circ w^i \\ ... \\ a_i^n \circ w^i \end{pmatrix}$

(6.4) $\quad a_j^i = a_{j\,s0}^i \otimes a_{j\,s1}^i$

**Summary of Results 6.3.** *Let passive transformation $g$ map the basis $\overline{\overline{e}}_1$ into the basis $\overline{\overline{e}}_2$:* (6.6) $\quad \overline{\overline{e}}_2 = g^*{}_* \overline{\overline{e}}_1$

*Let* (6.7) $\quad a_{k\,j}^{\phantom{k}i} = a_{k\,j\,s0}^{\phantom{k}i} \otimes a_{k\,j\,s1}^{\phantom{k}i}$ *be coordinates of linear map $a$ with respect to the basis $\overline{\overline{e}}_k$, $k = 1, 2$. Then*

(6.13) $\quad a_{2\,k\,t0}^{\phantom{2k}l} \otimes a_{2\,k\,t1}^{\phantom{2k}l} = (a_{1\,j\,s0}^{\phantom{1j}i} \otimes g_k^j a_{1\,j\,s1}^{\phantom{1j}i}) g^{-1*}{}_*{}_i^l$

*Equalities*

(6.1) $\quad f_2 = g^*{}_* f_1^*{}_* g^{-1*}{}_*$

*and*

(6.13) $\quad a_{2\,k\,t0}^{\phantom{2k}l} \otimes a_{2\,k\,t1}^{\phantom{2k}l} = (a_{1\,j\,s0}^{\phantom{1j}i} \otimes g_k^j a_{1\,j\,s1}^{\phantom{1j}i}) g^{-1*}{}_*{}_i^l$

*are similar. From the equality*

(6.14) $\begin{aligned} a \circ v &= \underline{a_{2\,k\,t0}^{\phantom{2k}l} v_2^k a_{2\,k\,t1}^{\phantom{2k}l} e_{2l}} = (a_{1\,j\,s0}^{\phantom{1j}i} v_2^k g_k^j a_{1\,j\,s1}^{\phantom{1j}i}) g^{-1*}{}_*{}_i^l g_l^p e_{1p} \\ &= \underline{(a_{1\,j\,s0}^{\phantom{1j}i} v_1^j a_{1\,j\,s1}^{\phantom{1j}i}) e_{1i}} \end{aligned}$

*it follows that image of linear map does not depend on choice of basis.*

*Let the set of tensors* (6.19) $\quad a_{k\,j_1...j_n\,s0}^{\phantom{k}i} \otimes ... \otimes a_{k\,j_1...j_n\,sn}^{\phantom{k}i}$ *be coordinates of polylinear map $a$ with respect to the basis $\overline{\overline{e}}_k$, $k = 1, 2$. Then*

(6.23) $\begin{aligned} &a_{2\,k_1...k_n\,t0}^{\phantom{2}l} \otimes a_{2\,k_1...k_n\,t1}^{\phantom{2}l} ... \otimes a_{2\,k_1...k_n\,tn}^{\phantom{2}l} \\ &= (a_{1\,j_1...j_n\,s0}^{\phantom{1}i} \otimes g_{k_1}^{j_1} a_{1\,j_1...j_n\,s1}^{\phantom{1}i} ... \otimes g_{k_n}^{j_n} a_{1\,j_1...j_n\,sn}^{\phantom{1}i}) g^{-1*}{}_*{}_i^l \end{aligned}$

*From the equality*

(6.24) $\begin{aligned} a \circ (v_1, ..., v_2) &= \underline{a_{2\,k_1...k_n\,t0}^{\phantom{2}l} v_{12}^{k_1} a_{2\,k_1...k_n\,t1}^{\phantom{2}l} ... v_{n2}^{k_n} a_{2\,k_1...k_n\,tn}^{\phantom{2}l} e_{2l}} \\ &= (a_{1\,j_1...j_n\,s0}^{\phantom{1}i} v_{12}^{k_1} g_{k_1}^{j_1} a_{1\,j_1...j_n\,s1}^{\phantom{1}i} ... v_{n2}^{k_n} g_{k_n}^{j_n} a_{1\,j_1...j_n\,sn}^{\phantom{1}i}) g^{-1*}{}_*{}_i^l g_l^p e_{1p} \\ &= \underline{(a_{1\,j_1...j_n\,s0}^{\phantom{1}i} v_{11}^{j_1} a_{1\,j_1...j_n\,s1}^{\phantom{1}i} ... v_{n1}^{j_n} a_{1\,j_1...j_n\,sn}^{\phantom{1}i}) e_{1i}} \end{aligned}$

it follows that image of polylinear map does not depend on choice of basis.

Consider skew-symmetric polylinear map

$$(6.26) \quad \begin{aligned} h^i \circ (u \otimes v) &= a^i_{jk\,s0} u^j a^i_{jk\,s1} v^k a^i_{jk\,s2} \\ &= (a^i_{jk\,s0} \otimes a^i_{jk\,s1} \otimes a^i_{jk\,s2}) \circ (u^j \otimes v^k) \end{aligned}$$

If we consider determinant like expression

$$(6.28) \quad \det{}^* \begin{pmatrix} u^j & v^j \\ u^k & v^k \end{pmatrix} = u^j \otimes v^k - v^j \otimes u^k$$

then we can represent polylinear map $h$ as

$$(6.29) \quad h^i \circ (u \otimes v) = \tfrac{1}{2}(a^i_{jk\,s0} \otimes a^i_{jk\,s1} \otimes a^i_{jk\,s2}) \circ \det{}^* \begin{pmatrix} u^j & v^j \\ u^k & v^k \end{pmatrix}$$

Let the set of tensors $(6.19) \quad a_k{}^i_{j_1 \ldots j_n\,s0} \otimes \ldots \otimes a_k{}^i_{j_1 \ldots j_n\,sn}$ be coordinates of polylinear map $h$ with respect to the basis $\overline{\overline{e}}_k$, $k = 1, 2$. Then

$$(6.37) \quad \begin{aligned} (a_2{}^i_{jk\,s0} \otimes_1 a_2{}^i_{jk\,s1} \otimes_2 a_2{}^i_{jk\,s2}) \circ \det{}^* \begin{pmatrix} 1 \otimes_1 \delta^j_q & 1 \otimes_2 \delta^j_r \\ 1 \otimes_1 \delta^k_q & 1 \otimes_2 \delta^k_r \end{pmatrix} \\ = ((a_1{}^i_{jk\,s0} \otimes_1 a_1{}^i_{jk\,s1} \otimes_2 a_1{}^i_{jk\,s2}) g^{-1*}{}^{*p}_i) \circ \det{}^* \begin{pmatrix} 1 \otimes_1 g^j_q & 1 \otimes_2 g^j_r \\ 1 \otimes_1 g^k_q & 1 \otimes_2 g^k_r \end{pmatrix} \end{aligned}$$

□

Let passive transformation $g$ map the basis $\overline{\overline{e}}_1$ into the basis $\overline{\overline{e}}_2$

$$(6.6) \quad \overline{\overline{e}}_2 = g^*{}_* \overline{\overline{e}}_1$$

Let

$$(6.7) \quad a_k{}^i_j = a_k{}^i_{j\,s0} \otimes a_k{}^i_{j\,s1}$$

be coordinates of linear map $a$ with respect to the basis $\overline{\overline{e}}_k$, $k = 1, 2$.

Let $v^j_k$, $w^j_k$ be coordinates of vectors $v$, $w$ with respect to the basis $\overline{\overline{e}}_k$, $k = 1, 2$. Then

$$(6.8) \quad w^i_k = a_k{}^i_{j\,s0} v^j_k a_k{}^i_{j\,s1}$$

According to the theorem 6.1,

$$(6.9) \quad v_1 = v_2{}^*{}_* g$$

$$(6.10) \quad w_1 = w_2{}^*{}_* g$$

Equalities

$$(6.11) \quad w^l_2 g^i_l = a_1{}^i_{j\,s0} v^k_2 g^j_k a_1{}^i_{j\,s1}$$

$$(6.12) \quad w^l_2 = (a_1{}^i_{j\,s0} v^k_2 g^j_k a_1{}^i_{j\,s1}) g^{-1*}{}^{*l}_i$$

follow from equalities

$(6.9) \quad v_1 = v_2{}^*{}_* g \qquad (6.10) \quad w_1 = w_2{}^*{}_* g$

(6.8)  $w_1^i = a_{1\,j}^{\,\,i}{}_{s0} v_1^j a_{1\,j}^{\,\,i}{}_{s1}$   (k=1).

The equality

(6.13) $$a_{2\,k}^{\,\,l}{}_{t0} \otimes a_{2\,k}^{\,\,l}{}_{t1} = (a_{1\,j}^{\,\,i}{}_{s0} \otimes g_k^j a_{1\,j}^{\,\,i}{}_{s1}) g^{-1*}{}^{*}{}_i^l$$

follows from equalities

(6.8)  $w_2^i = a_{2\,j}^{\,\,i}{}_{s0} v_2^j a_{2\,j}^{\,\,i}{}_{s1}$   (k=2),

(6.12)  $w_2^l = (a_{1\,j}^{\,\,i}{}_{s0} v_2^k g_k^j a_{1\,j}^{\,\,i}{}_{s1}) g^{-1*}{}^{*}{}_i^l$

The equality

(6.14) $$\begin{aligned} a \circ v &= a_{2\,k}^{\,\,l}{}_{t0} v_2^k a_{2\,k}^{\,\,l}{}_{t1} e_{2l} = (a_{1\,j}^{\,\,i}{}_{s0} v_2^k g_k^j a_{1\,j}^{\,\,i}{}_{s1}) g^{-1*}{}^{*}{}_i^l g_l^p e_{1p} \\ &= (a_{1\,j}^{\,\,i}{}_{s0} v_1^j a_{1\,j}^{\,\,i}{}_{s1}) e_{1i} \end{aligned}$$

follows from equalities

(6.9)  $v_1 = v_2{}^{*}{}_{*} g$

(6.13)  $a_{2\,k}^{\,\,l}{}_{t0} \otimes a_{2\,k}^{\,\,l}{}_{t1} = (a_{1\,j}^{\,\,i}{}_{s0} \otimes g_k^j a_{1\,j}^{\,\,i}{}_{s1}) g^{-1*}{}^{*}{}_i^l$

From the equality (6.14) it follows that image of linear map does not depend on choice of basis.

It is easy to see the similarity of equalities

(6.1)  $f_2 = g^*{}_{*} f_1{}^{*}{}_{*} g^{-1*}{}^{*}$

and

(6.13)  $a_{2\,k}^{\,\,l}{}_{t0} \otimes a_{2\,k}^{\,\,l}{}_{t1} = (a_{1\,j}^{\,\,i}{}_{s0} \otimes g_k^j a_{1\,j}^{\,\,i}{}_{s1}) g^{-1*}{}^{*}{}_i^l$

If $D$-algebra $A$ is commutative, then these equalities are the same.

**Theorem 6.4.** *Polylinear map*

(6.15) $$a : V^n \to V$$

*is geometric object.*

PROOF. Without loss of generality, we consider proof of the theorem for left $A$-vector space of columns.

Let passive transformation $g$ map the basis $\overline{\overline{e}}_1$ into the basis $\overline{\overline{e}}_2$

(6.16) $$\overline{\overline{e}}_2 = g^*{}_{*} \overline{\overline{e}}_1$$

Let $v_{ik}^j$, $w_k^j$ be coordinates of vectors $v_i$, $i = 1, ..., n$, $w$ with respect to the basis $\overline{\overline{e}}_k$, $k = 1, 2$. According to the theorem 6.1,

(6.17) $$v_{i1} = v_{i2}{}^{*}{}_{*} g$$

(6.18) $$w_1 = w_2{}^{*}{}_{*} g$$

Let the set of tensors

(6.19) $$a_{k\,j_1...j_n}^{\,\,i}{}_{s0} \otimes ... \otimes a_{k\,j_1...j_n}^{\,\,i}{}_{sn}$$

be coordinates of polylinear map $a$ with respect to the basis $\overline{\overline{e}}_k$. Then[6.1]

(6.20) $$w_k^i = a_{kj_1\ldots j_n s0}^i v_{1k}^{j_1} a_{kj_1\ldots j_n s1}^i \ldots v_{nk}^{j_n} a_{kj_1\ldots j_n sn}^i$$

Equalities

(6.21) $$w_2^l g_l^i = a_{1j_1\ldots j_n s0}^i v_{12}^{k_1} g_{k_1}^{j_1} a_{1j_1\ldots j_n s1}^i \ldots v_{n2}^{k_n} g_{k_n}^{j_n} a_{kj_1\ldots j_n sn}^i$$

(6.22) $$w_2^l = (a_{1j_1\ldots j_n s0}^i v_{12}^{k_1} g_{k_1}^{j_1} a_{1j_1\ldots j_n s1}^i \ldots v_{n2}^{k_n} g_{k_n}^{j_n} a_{kj_1\ldots j_n sn}^i) g^{-1*}{}_*{}^l{}_i$$

follow from equalities

$\boxed{(6.17) \quad v_{i1} = v_{i2}{}^*{}_* g}$ $\boxed{(6.18) \quad w_1 = w_2{}^*{}_* g}$

$\boxed{(6.20) \quad w_1^i = a_{1j_1\ldots j_n s0}^i v_{11}^{j_1} a_{1j_1\ldots j_n s1}^i \ldots v_{n1}^{j_n} a_{1j_1\ldots j_n sn}^i}$ (k=1).

The equality

(6.23) $$a_{2k_1\ldots k_n t0}^l \otimes a_{2k_1\ldots k_n t1}^l \ldots \otimes a_{2k_1\ldots k_n tn}^l$$
$$= (a_{1j_1\ldots j_n s0}^i \otimes g_{k_1}^{j_1} a_{1j_1\ldots j_n s1}^i \ldots \otimes g_{k_n}^{j_n} a_{1j_1\ldots j_n sn}^i) g^{-1*}{}_*{}^l{}_i$$

follows from equalities

$\boxed{(6.20) \quad w_2^i = a_{2j_1\ldots j_n s0}^i v_{12}^{j_1} a_{2j_1\ldots j_n s1}^i \ldots v_{n2}^{j_n} a_{2j_1\ldots j_n sn}^i}$ (k=2),

$\boxed{(6.22) \quad w_2^l = (a_{1j_1\ldots j_n s0}^i v_{12}^{k_1} g_{k_1}^{j_1} a_{1j_1\ldots j_n s1}^i \ldots v_{n2}^{k_n} g_{k_n}^{j_n} a_{kj_1\ldots j_n sn}^i) g^{-1*}{}_*{}^l{}_i}$

The equality

(6.24)
$$a \circ (v_1, \ldots, v_2) = \underline{a_{2k_1\ldots k_n t0}^l v_{12}^{k_1} a_{2k_1\ldots k_n t1}^l \ldots v_{n2}^{k_n} a_{2k_1\ldots k_n tn}^l e_{2l}}$$
$$= (a_{1j_1\ldots j_n s0}^i v_{12}^{k_1} g_{k_1}^{j_1} a_{1j_1\ldots j_n s1}^i \ldots v_{n2}^{k_n} g_{k_n}^{j_n} a_{1j_1\ldots j_n sn}^i) g^{-1*}{}_*{}^l{}_i g_l^p e_{1p}$$
$$= \underline{(a_{1j_1\ldots j_n s0}^i v_{11}^{j_1} a_{1j_1\ldots j_n s1}^i \ldots v_{n1}^{j_n} a_{1j_1\ldots j_n sn}^i) e_{1i}}$$

follows from equalities

$\boxed{(6.17) \quad v_{i1} = v_{i2}{}^*{}_* g}$

$\boxed{\begin{aligned}(6.23) \quad & a_{2k_1\ldots k_n t0}^l \otimes a_{2k_1\ldots k_n t1}^l \ldots \otimes a_{2k_1\ldots k_n tn}^l \\ & = (a_{1j_1\ldots j_n s0}^i \otimes g_{k_1}^{j_1} a_{1j_1\ldots j_n s1}^i \ldots \otimes g_{k_n}^{j_n} a_{1j_1\ldots j_n sn}^i) g^{-1*}{}_*{}^l{}_i\end{aligned}}$

From the equality (6.24) it follows that image of polylinear map does not depend on choice of basis. □

We consider transformation skew-symmetric polylinear map using the bilinear map as an example. Let

(6.25) $$w = h \circ (u \otimes v) = \frac{1}{2} h \circ (u \otimes v - v \otimes u)$$

where the map $h$ has the form

(6.26) $$h^i \circ (u \otimes v) = a_{jk\,s0}^i u^j a_{jk\,s1}^i v^k a_{jk\,s2}^i$$
$$= (a_{jk\,s0}^i \otimes a_{jk\,s1}^i \otimes a_{jk\,s2}^i) \circ (u^j \otimes v^k)$$

---

[6.1] In this proof, there is no need to consider posibility to change an order of variables in every term since we are interested in single term. If the order of variables in this term changes, then we can change this order before we start calculations in the proof.

The equality

$$(6.27) \quad w^i = h^i \circ (u \otimes v) = \frac{1}{2}(a^i_{jk\,s0} \otimes_1 a^i_{jk\,s1} \otimes_2 a^i_{jk\,s2}) \circ (u^j \otimes v^k - v^j \otimes u^k)$$

follows from equalities (6.25), (6.26). If we consider determinant like expression

$$(6.28) \quad \det{}^* \begin{pmatrix} u^j & v^j \\ u^k & v^k \end{pmatrix} = u^j \otimes v^k - v^j \otimes u^k$$

then the equality (6.27) gets the form

$$(6.29) \quad h^i \circ (u \otimes v) = \frac{1}{2}(a^i_{jk\,s0} \otimes a^i_{jk\,s1} \otimes a^i_{jk\,s2}) \circ \det{}^* \begin{pmatrix} u^j & v^j \\ u^k & v^k \end{pmatrix}$$

Let passive transformation $g$ map the basis $\overline{\overline{e}}_1$ into the basis $\overline{\overline{e}}_2$

$$(6.30) \quad \overline{\overline{e}}_2 = g^*{}_*\overline{\overline{e}}_1$$

The equality

$$(6.31) \quad \begin{aligned} w_2^p g_p^i &= \frac{1}{2}(a_1{}^i_{jk\,s0} \otimes_1 a_1{}^i_{jk\,s1} \otimes_2 a_1{}^i_{jk\,s2}) \\ &\quad \circ ((u_2^q g_q^j) \otimes (v_2^r g_r^k) - (v_2^r g_r^j) \otimes (u_2^q g_q^k)) \\ &= \frac{1}{2}(a_1{}^i_{jk\,s0} \otimes_1 a_1{}^i_{jk\,s1} \otimes_2 a_1{}^i_{jk\,s2}) \\ &\quad \circ ((1 \otimes_1 g_q^j) \otimes (1 \otimes_2 g_r^k) - (1 \otimes_2 g_r^j) \otimes (1 \otimes_1 g_q^k)) \circ (u_2^q \otimes v_2^r) \end{aligned}$$

follows from equalities

$$\boxed{(6.27) \quad w^i = h^i \circ (u \otimes v) = \tfrac{1}{2}(a^i_{jk\,s0} \otimes_1 a^i_{jk\,s1} \otimes_2 a^i_{jk\,s2}) \circ (u^j \otimes v^k - v^j \otimes u^k)} \quad \boxed{(6.18) \quad w_1 = w_2{}^*{}_*g}$$

$$(6.32) \quad u_1 = u_2{}^*{}_*g$$

$$(6.33) \quad u_1 = u_2{}^*{}_*g$$

The equality

$$(6.34) \quad \begin{aligned} w_2^p &= \frac{1}{2}((a_1{}^i_{jk\,s0} \otimes_1 a_1{}^i_{jk\,s1} \otimes_2 a_1{}^i_{jk\,s2})g^{-1}{}_*{}^*{}_i^p) \\ &\quad \circ ((1 \otimes_1 g_q^j) \otimes (1 \otimes_2 g_r^k) - (1 \otimes_2 g_r^j) \otimes (1 \otimes_1 g_q^k)) \circ (u_2^q \otimes v_2^r) \end{aligned}$$

follows from the equality (6.31). The equality

$$(6.35) \quad \begin{aligned} w_2^p &= \frac{1}{2}(a_2{}^i_{jk\,s0} \otimes_1 a_2{}^i_{jk\,s1} \otimes_2 a_2{}^i_{jk\,s2}) \\ &\quad \circ ((1 \otimes_1 \delta_q^j) \otimes (1 \otimes_2 \delta_r^k) - (1 \otimes_2 \delta_r^j) \otimes (1 \otimes_1 \delta_q^k)) \circ (u_2^q \otimes v_2^r) \end{aligned}$$

follows from the equality

$$\boxed{(6.27) \quad w^i = h^i \circ (u \otimes v) = \tfrac{1}{2}(a^i_{jk\,s0} \otimes_1 a^i_{jk\,s1} \otimes_2 a^i_{jk\,s2}) \circ (u^j \otimes v^k - v^j \otimes u^k)}$$

The equality

$$(6.36) \quad \begin{aligned} &(a_2{}^i_{jk\,s0} \otimes_1 a_2{}^i_{jk\,s1} \otimes_2 a_2{}^i_{jk\,s2}) \\ &\circ ((1 \otimes_1 \delta_q^j) \otimes (1 \otimes_2 \delta_r^k) - (1 \otimes_2 \delta_r^j) \otimes (1 \otimes_1 \delta_q^k)) \\ &= ((a_1{}^i_{jk\,s0} \otimes_1 a_1{}^i_{jk\,s1} \otimes_2 a_1{}^i_{jk\,s2})g^{-1}{}_*{}^*{}_i^p) \\ &\circ ((1 \otimes_1 g_q^j) \otimes (1 \otimes_2 g_r^k) - (1 \otimes_2 g_r^j) \otimes (1 \otimes_1 g_q^k)) \end{aligned}$$

follows from the equality (6.34), (6.35). The equality

$$
\begin{aligned}
(6.37) \quad & (a_2{}^i_{jk}{}_{s0} \otimes_1 a_2{}^i_{jk}{}_{s1} \otimes_2 a_2{}^i_{jk}{}_{s2}) \circ \det{}^* \begin{pmatrix} 1 \otimes_1 \delta^j_q & 1 \otimes_2 \delta^j_r \\ 1 \otimes_1 \delta^k_q & 1 \otimes_2 \delta^k_r \end{pmatrix} \\
& = ((a_1{}^i_{jk}{}_{s0} \otimes_1 a_1{}^i_{jk}{}_{s1} \otimes_2 a_1{}^i_{jk}{}_{s2}) g^{-1}{}_*{}^{*\,p}_{\phantom{*}i}) \circ \det{}^* \begin{pmatrix} 1 \otimes_1 g^j_q & 1 \otimes_2 g^j_r \\ 1 \otimes_1 g^k_q & 1 \otimes_2 g^k_r \end{pmatrix}
\end{aligned}
$$

follows from the equality (6.28), (6.36).

## 8. Index

# 9. Special Symbols and Notations

# Ковариантность в некомутативной алгебре

Александр Клейн

Аннотация. Рассмотрим векторное пространство над некоммутативной алгеброй с делением. Множество автоморфизмов этого векторного пространства является группой $GL$. Группа $GL$ действует на множестве базисов векторного пространства (многообразие базисов) однотранзитивно и порождает активное представление. Парное представление на многообразии базисов называется пассивным представлением. Не существует автоморфизма, ассоциированного с пассивным преобразованием. Однако пассивное преобразование порождает преобразование координат вектора относительно базиса. Если мы рассмотрим гомоморфизм векторного пространства $V$ в векторное пространство $W$, то мы можем изучить как пассивное преобразование в векторном пространстве $V$ порождает преобразование координат вектора в векторном пространстве $W$. Вектор в векторном пространстве $W$ называется геометрическим объектом в векторном пространстве $V$. Принцип ковариантности утверждает, что геометрический объект не зависит от выбора базиса. Я рассмотрел преобразование координат вектора и полилинейного отображения.

## Содержание

## 1. Предисловие

Рассмотрим резиновую трубу, направленную вертикально. Воздух, который движется внутри трубы, меняет поверхность трубы таким образом, что горизонтальное сечение имеет вид

$$\frac{x^2}{a^2} + \frac{y^2}{b^2} = 1 \tag{1.1}$$

Aleks_Kleyn@MailAPS.org.
  http://AleksKleyn.dyndns-home.com:4080/   http://arxiv.org/a/kleyn_a_1.
  http://AleksKleyn.blogspot.com/.

где параметры $a$, $b$ являются функцией координаты $z$ и времени $t$ и являются решением системы дифференциальных уравнений

$$
\begin{aligned}
\frac{\partial a}{\partial z} &= a_1 x + b_1 y & \frac{\partial b}{\partial z} &= c_1 x + d_1 y \\
\frac{\partial a}{\partial t} &= a_2 x + b_2 y & \frac{\partial b}{\partial t} &= c_2 x + d_2 y
\end{aligned}
\tag{1.2}
$$

Однако расстояние может измеряться в метрах или сантиметрах, время может измеряться в минутах или секундах. Мы можем изменить направление осей $x$, $y$ в плоскости, перпендикулярной оси $z$. Эти преобразования координат отразятся в записи уравнения и его решения. Однако решение задачи не зависит от того, какой вид это решение имеет в координатном представлении.

Рассмотрим более простой пример. Кривая второго порядка на плоскости имеет очень сложный вид. Однако мы всегда можем выбрать координаты, относительно которых кривая принимает канонический вид и является элипсом, гиперболой, параболой или парой прямых.

Симметрия, состоящая в независимости решения от выбора координат, называется ковариантностью.

Тензорная алгебра является выражением ковариантности в векторном пространстве над полем. Наша задача рассмотреть теорию ковариантности в векторном пространстве над $D$-алгеброй. Хотя эта теория сложнее чем теория ковариантности в векторном пространстве над полем, обе теории имеют общие черты.

## 2. Произведение матриц в некоммутативной алгебре

Прежде чем мы начнём рассматривать векторные пространства над некоммутативной алгеброй $A$ с делением, я напомню, что существует две операции произведения матриц с элементами из алгебры $A$.

**Определение 2.1.** *Пусть число столбцов матрицы $a$ равно числу строк матрицы $b$. $_*^*$-произведение матриц $a$ и $b$ имеет вид*

$$a_*^*b = \left(a^i_k b^k_j\right) \tag{2.1}$$

$$(a_*^*b)^i_j = a^i_k b^k_j \tag{2.2}$$

$$
\begin{pmatrix} a^1_1 & \ldots & a^1_p \\ \ldots & \ldots & \ldots \\ a^n_1 & \ldots & a^n_p \end{pmatrix} {}_*^* \begin{pmatrix} b^1_1 & \ldots & b^1_m \\ \ldots & \ldots & \ldots \\ b^p_1 & \ldots & b^p_m \end{pmatrix} = \begin{pmatrix} a^1_k b^k_1 & \ldots & a^1_k b^k_m \\ \ldots & \ldots & \ldots \\ a^n_k b^k_1 & \ldots & a^n_k b^k_m \end{pmatrix}
$$
$$
= \begin{pmatrix} (a_*^*b)^1_1 & \ldots & (a_*^*b)^1_m \\ \ldots & \ldots & \ldots \\ (a_*^*b)^n_1 & \ldots & (a_*^*b)^n_m \end{pmatrix}
\tag{2.3}
$$

*$_*^*$-произведение может быть выражено как произведение строк матрицы $a$ и столбцов матрицы $b$.* □

**Определение 2.2.** *Пусть число строк матрицы $a$ равно числу столбцов матрицы $b$. $^*_*$-**произведение** матриц $a$ и $b$ имеет вид*

(2.4) $$a^*{}_*b = \left(a_i^k b_k^j\right)$$

(2.5) $$(a^*{}_*b)_j^i = a_i^k b_k^j$$

(2.6) $$\begin{pmatrix} a_1^1 & ... & a_m^1 \\ ... & ... & ... \\ a_1^p & ... & a_m^p \end{pmatrix} *_* \begin{pmatrix} b_1^1 & ... & b_p^1 \\ ... & ... & ... \\ b_1^n & ... & b_p^n \end{pmatrix} = \begin{pmatrix} a_1^k b_k^1 & ... & a_m^k b_k^1 \\ ... & ... & ... \\ a_1^k b_k^n & ... & a_m^k b_k^n \end{pmatrix}$$
$$= \begin{pmatrix} (a^*{}_*b)_1^1 & ... & (a^*{}_*b)_m^1 \\ ... & ... & ... \\ (a^*{}_*b)_1^n & ... & (a^*{}_*b)_m^n \end{pmatrix}$$

*$^*_*$-произведение может быть выражено как произведение столбцов матрицы $a$ и строк матрицы $b$.* □

## 3. Гомоморфизм левого $A$-векторного пространства

Пусть $A$ - некоммутативная $D$-алгебра с делением. Тогда мы можем рассматривать левое или правое $A$-векторное пространство. Мы можем ограничиться рассмотрением левого $A$-векторного пространства, так как утверждения для правого $A$-векторного пространства аналогичны.

Мы начнём с рассмотрения множества гомомрфизмов левого $A$-векторного пространства. Гомоморфизм - это отображение, которое сохраняет структуру алгебры. В частности, верна теорема 3.2.

**Определение 3.1.** *Пусть*

(3.1) 
$$A \xrightarrow{g_{1.23}} A \xrightarrow{g_{1.34}} V_1 \qquad \begin{aligned} g_{12}(d) &: a \to d\,a \\ g_{23}(v) &: w \to C(w,v) \\ C &\in \mathcal{L}(A^2 \to A) \\ g_{3,4}(a) &: v \to a\,v \\ g_{1,4}(d) &: v \to d\,v \end{aligned}$$

с участием $g_{1.12}$, $g_{1.12}$, $g_{1.14}$, $D$.

*диаграмма представлений, описывающая левый $A$-модуль $V_1$. Пусть*

(3.2)
$$A \xrightarrow{g_{2.23}} A \xrightarrow{g_{2.34}} V_2 \qquad \begin{aligned} g_{12}(d) &: a \to d\,a \\ g_{23}(v) &: w \to C(w,v) \\ C &\in \mathcal{L}(A^2 \to A) \\ g_{3,4}(a) &: v \to a\,v \\ g_{1,4}(d) &: v \to d\,v \end{aligned}$$

с участием $g_{2.12}$, $g_{2.12}$, $g_{2.14}$, $D$.

*диаграмма представлений, описывающая левый $A$-модуль $V_2$. Морфизм*

$$(3.3) \qquad f : V_1 \to V_2$$

*диаграммы представлений* (3.1) *в диаграмму представлений* (3.2) *называется* **гомоморфизмом** *левого $A$-модуля $V_1$ в левый $A$-модуль $V_2$. Обозначим* $\text{Hom}(D; A; *V_1 \to *V_2)$ *множество гомоморфизмов левого $A$-модуля $V_1$ в левый $A$-модуль $V_2$.* □

Мы будем пользоваться записью

$$f \circ a = f(a)$$

для образа гомоморфизма $f$.

**Теорема 3.2.** *Гомоморфизм*

$$(3.3) \quad f : V_1 \to V_2$$

*левого $A$-модуля $V_1$ в левый $A$-модуль $V_2$ удовлетворяет следующим равенствам*

$$(3.4) \qquad f \circ (u + v) = f \circ u + f \circ v$$

$$(3.5) \qquad f \circ (av) = a(f \circ v)$$

$$a \in A \quad u, v \in V_1$$

Доказательство. Теорема является следствием теоремы [3]-10.3.2. □

**Определение 3.3.** *Гомоморфизм*

$$f : V \to W$$

*называется*[3.1]**изоморфизмом** *между левыми $A$-векторными пространствами $V$ и $W$, если соответствие $f^{-1}$ является гомоморфизмом. Гомоморфизм*

$$f : V \to V$$

*источником и целью которого является одно и тоже левое $A$-векторное пространство, называется* **эндоморфизмом**. *Эндоморфизм*

$$f : V \to V$$

*левого $A$-векторного пространства $V$ называется* **автоморфизмом**, *если соответствие $f^{-1}$ является эндоморфизмом.* □

**Теорема 3.4.** *Множество $GL(V)$ автоморфизмов левого $A$-векторного пространства $V$ является группой.*

Хотя теорема 3.2 верна и может служить определением гомоморфизма, эта теорема не даёт полной картины.

Допустим мы изучаем физическое явление. Мы можем дать качественную картину этого явления. Но если мы хотим провести эксперимент, мы должны выбрать набор измерительных инструментов и рассчитать ожидаемые измерения. Набор измерительных инструментов является базисом в пространстве

---

[3.1] Я следую определению на странице [4]-63.

измерений, а ожидаемые измерения являются координатами эксперимента в пространстве измерений.

Мы распространим это замечание на теорему 3.2. Для того, чтобы понять структуру гомоморфизма и ответить на вопрос как велико множество гомоморфизмов, мы должны выбрать базис $A$-векторного пространства и рассмотреть изменение координат вектора относительно выбраного базиса.

Существуют различные соглашения о нумерации векторов базиса и соответствующие соглашения о нумерации координат вектора. Например, в теореме 3.5 мы рассматриваем координаты вектора в левом $A$-векторном пространстве столбцов.

**Теорема 3.5.** *Если мы запишем векторы базиса $\overline{\overline{e}}$ в виде строки матрицы*

$$(3.6) \qquad e = \begin{pmatrix} e_1 & ... & e_n \end{pmatrix}$$

*и координаты вектора $\overline{w} = w^i e_i$ относительно базиса $\overline{\overline{e}}$ в виде столбца матрицы*

$$(3.7) \qquad w = \begin{pmatrix} w^1 \\ ... \\ w^n \end{pmatrix}$$

*то мы можем представить вектор $\overline{w}$ в виде $^*{}_*$-произведения матриц*

$$(3.8) \qquad \overline{w} = w^*{}_* e = \begin{pmatrix} w^1 \\ ... \\ w^n \end{pmatrix} {}^*{}_* \begin{pmatrix} e_1 & ... & e_n \end{pmatrix} = w^i e_i$$

Если базис $\overline{\overline{e}}$ задан, то

$$(3.9) \qquad V_* = \{v : \overline{v} = v^*{}_* e, \overline{v} \in V\}$$

является множеством координат векторов $\overline{v} \in V$. Множество $V_*$ является левым $A$-векторным пространством и изоморфно левому $A$-векторному пространству $V$. Следовательно, множество $V_*$ не зависит от выбора базиса $\overline{\overline{e}}$.

Если левое $A$-векторное пространство $V$ имеет размерность $n$, то левое $A$-векторное пространство $V_*$ изоморфно прямой сумме $n$ копий $D$-алгебры $A$. Поэтому мы в этом случае положим

$$(3.10) \qquad V_* = n^*{}_* A$$

В теореме 3.6 мы рассматриваем координаты вектора в левом $A$-векторном пространстве строк.

**Теорема 3.6.** *Если мы запишем векторы базиса $\overline{\overline{e}}$ в виде столбца матрицы*

$$(3.11) \qquad e = \begin{pmatrix} e^1 \\ ... \\ e^n \end{pmatrix}$$

и координаты вектора $\overline{w} = w_i e^i$ относительно базиса $\overline{\overline{e}}$ в виде строки матрицы

$$(3.12) \qquad w = \begin{pmatrix} w_1 & \dots & w_n \end{pmatrix}$$

то мы можем представить вектор $\overline{w}$ в виде $_*^*$-произведения матриц

$$(3.13) \qquad \overline{w} = w_*{}^* e = \begin{pmatrix} w_1 & \dots & w_n \end{pmatrix} {}_*^* \begin{pmatrix} e^1 \\ \dots \\ e^n \end{pmatrix} = w_i e^i$$

Мы рассмотрели основные форматы представления координат вектора левого $A$-векторного пространства. Очевидно, что возможны и другие форматы представления координат вектора. Например, мы можем рассматривать множество $n \times m$ матриц как левое $A$-векторное пространство. Аналогичным образом мы можем рассмотреть основные форматы представления координат вектора правого $A$-векторного пространства.

Утверждения для различных форм представления координат вектора похожи. Поэтому мы сосредоточим наше внимание на левом $A$-векторном пространстве столбцов.

**Теорема 3.7.** *Гомоморфизм*[3.2]

$$(3.14) \qquad \overline{f} : V_1 \to V_2$$

*левого $A$-векторного пространства $V_1$ в левое $A$-векторное пространство $V_2$ имеет представление*

$$(3.15) \qquad w = v^*{}_* f$$

$$(3.16) \qquad \overline{f} \circ (v^i e_{V_1 i}) = v^i f_i^k e_{V_2 k}$$

$$(3.17) \qquad \overline{f} \circ (v^*{}_* e_1) = v^*{}_* f^*{}_* e_2$$

*относительно выбранных базисов. Здесь*

- $v$ - координатная матрица $V_1$-числа $\overline{v}$ относительно базиса $\overline{\overline{e}}_{V_1}$

$$(3.18) \qquad \overline{v} = v^*{}_* e_{V_1}$$

- $w$ - координатная матрица $V_2$-числа

$$(3.19) \qquad \overline{w} = \overline{f} \circ \overline{v}$$

**Теорема 3.8.** *Гомоморфизм*[3.3]

$$(3.21) \qquad \overline{f} : V_1 \to V_2$$

*левого $A$-векторного пространства $V_1$ в левое $A$-векторное пространство $V_2$ имеет представление*

$$(3.22) \qquad w = v_*{}^* f$$

$$(3.23) \qquad \overline{f} \circ (v_i e_{V_1}^i) = v_i f_k^i e_{V_2}^k$$

$$(3.24) \qquad \overline{f} \circ (v_*{}^* e_1) = v_*{}^* f_*{}^* e_2$$

*относительно выбранных базисов. Здесь*

- $v$ - координатная матрица $V_1$-числа $\overline{v}$ относительно базиса $\overline{\overline{e}}_{V_1}$

$$(3.25) \qquad \overline{v} = v_*{}^* e_{V_1}$$

- $w$ - координатная матрица $V_2$-числа

$$(3.26) \qquad \overline{w} = \overline{f} \circ \overline{v}$$

---

[3.2] В теоремах 3.7, 3.9, мы опираемся на следующее соглашение. Пусть множество векторов $\overline{\overline{e}}_1 = (e_{1i}, i \in I)$ является базисом левого $A$-векторного пространства $V_1$. Пусть множество векторов $\overline{\overline{e}}_2 = (e_{2j}, j \in J)$ является базисом левого $A$-векторного пространства $V_2$.

[3.3] В теоремах 3.8, 3.10, мы опираемся на следующее соглашение. Пусть множество векторов $\overline{\overline{e}}_1 = (e_1^i, i \in I)$ является базисом левого $A$-векторного пространства $V_1$. Пусть множество векторов $\overline{\overline{e}}_2 = (e_2^j, j \in J)$ является базисом левого $A$-векторного пространства $V_2$.

| | |
|---|---|
| *относительно базиса* $\overline{\overline{e}}_{V_2}$<br>(3.20) $\qquad \overline{w} = w^*{}_* e_{V_2}$<br>- $f$ - *координатная матрица множества* $V_2$*-чисел* $(\overline{f} \circ e_{V_1 i}, i \in I)$ *относительно базиса* $\overline{\overline{e}}_{V_2}$.<br><br>*Матрица* $f$ *определена однозначно и называется* **матрицей гомоморфизма** $\overline{f}$ *относительно базисов* $\overline{\overline{e}}_1$, $\overline{\overline{e}}_2$. | *относительно базиса* $\overline{\overline{e}}_{V_2}$<br>(3.27) $\qquad \overline{w} = w_*{}^* e_{V_2}$<br>- $f$ - *координатная матрица множества* $V_2$*-чисел* $(\overline{f} \circ e_{V_1}^i, i \in I)$ *относительно базиса* $\overline{\overline{e}}_{V_2}$.<br><br>*Матрица* $f$ *определена однозначно и называется* **матрицей гомоморфизма** $\overline{f}$ *относительно базисов* $\overline{\overline{e}}_1$, $\overline{\overline{e}}_2$. |

Доказательство. Теорема является следствием теоремы [3]-10.3.3. □

Обратная теорема также верна.

Доказательство. Теорема является следствием теоремы [3]-10.3.4. □

| | |
|---|---|
| **Теорема 3.9.** *Пусть*<br>$$f = (f_j^i, i \in I, j \in J)$$<br>*матрица* $A$*-чисел. Отображение*[3.2]<br>(3.14) $\quad \overline{f} : V_1 \to V_2$<br>*определённое равенством*<br>(3.17) $\quad \overline{f} \circ (v^*{}_* e_1) = v^*{}_* f^*{}_* e_2$<br>*является гомоморфизмом левого* $A$*-векторного пространства столбцов. Гомоморфизм* (3.14), *который имеет данную матрицу* $f$, *определён однозначно.* | **Теорема 3.10.** *Пусть*<br>$$f = (f_i^j, i \in I, j \in J)$$<br>*матрица* $A$*-чисел. Отображение*[3.3]<br>(3.21) $\quad \overline{f} : V_1 \to V_2$<br>*определённое равенством*<br>(3.24) $\quad \overline{f} \circ (v_*{}^* e_1) = v_*{}^* f_*{}^* e_2$<br>*является гомоморфизмом левого* $A$*-векторного пространства строк. Гомоморфизм* (3.21), *который имеет данную матрицу* $f$, *определён однозначно.* |

Доказательство. Теорема является следствием теоремы [3]-10.3.5. □

Доказательство. Теорема является следствием теоремы [3]-10.3.6. □

Следовательно, если выбран базис $\overline{\overline{e}} = (e_1, ..., e_n)$ левого $A$-векторного пространства $V$ столбцов, то мы можем отождествить множество эндоморфизмов левого $A$-векторного пространства $V$ и множество $n \times n$ матриц. Эта связь между эндоморфизмами и матрицами будет видна глубже, если мы рассмотрим автоморфизмы.

| | |
|---|---|
| **Теорема 3.11.** *Пусть* $V$ - *левое* $A$*-векторное пространство столбцов и* $\overline{\overline{e}}$ - *базис левого* $A$*-векторного пространства* $V$. *Любой автоморфизм* $\overline{f}$ *левого* $A$*-векторного пространства* $V$ *имеет вид*<br>(3.28) $\qquad v' = v^*{}_* f$<br>*где* $f$ - $^*{}_*$-*невырожденная матрица. Матрицы автоморфизмов левого* $A$- | **Теорема 3.12.** *Пусть* $V$ - *левое* $A$*-векторное пространство строк и* $\overline{\overline{e}}$ - *базис левого* $A$*-векторного пространства* $V$. *Любой автоморфизм* $\overline{f}$ *левого* $A$*-векторного пространства* $V$ *имеет вид*<br>(3.30) $\qquad v' = v_*{}^* f$<br>*где* $f$ - $_*{}^*$-*невырожденная матрица. Матрицы автоморфизмов левого* $A$- |

*векторного пространства столбцов порождают группу $GL(V_*)$, изоморфную группе $GL(V)$. Автоморфизмы левого $A$-векторного пространства столбцов порождают правостороннее линейное эффективное представление*

(3.29) $\qquad GL(V_*) \xrightarrow{\phantom{xx}*\phantom{xx}} V_*$

*группы $GL(V_*)$ в левом $A$-векторном пространстве $V_*$.*

*векторного пространства строк порождают группу $GL(V_*)$, изоморфную группе $GL(V)$. Автоморфизмы левого $A$-векторного пространства строк порождают правостороннее линейное эффективное представление*

(3.31) $\qquad GL(V_*) \xrightarrow{\phantom{xx}*\phantom{xx}} V_*$

*группы $GL(V_*)$ в левом $A$-векторном пространстве $V_*$.*

Доказательство. Теорема является следствием теоремы [3]-12.1.9, [3]-12.1.11, [3]-12.2.1. $\square$

Доказательство. Теорема является следствием теоремы [3]-12.1.10, [3]-12.1.12, [3]-12.2.2. $\square$

Теоремы 3.14, 3.16 показывают связь между $_**$-невырожденными матрицами, базисами и автоморфизмами левого $A$-векторного пространства строк.

**Теорема 3.13.** *Пусть $V$ - левое $A$-векторное пространство столбцов. Координатная матрица базиса $\overline{\overline{g}}$ относительно базиса $\overline{\overline{e}}$ левого $A$-векторного пространства $V$ является $_**$-невырожденной матрицей.*

**Теорема 3.14.** *Пусть $V$ - левое $A$-векторное пространство строк. Координатная матрица базиса $\overline{\overline{g}}$ относительно базиса $\overline{\overline{e}}$ левого $A$-векторного пространства $V$ является $_**$-невырожденной матрицей.*

Доказательство. Теорема является следствием теоремы [3]-12.1.7. $\square$

Доказательство. Теорема является следствием теоремы [3]-12.1.8. $\square$

**Теорема 3.15.** *Автоморфизм $a$, действуя на каждый вектор базиса в левом $A$-векторном пространстве столбцов, отображает базис в другой базис.*

**Теорема 3.16.** *Автоморфизм $a$, действуя на каждый вектор базиса в левом $A$-векторном пространстве строк, отображает базис в другой базис.*

Доказательство. Теорема является следствием теоремы [3]-12.2.3. $\square$

Доказательство. Теорема является следствием теоремы [3]-12.2.4. $\square$

Таким образом, мы можем распространить правостороннее линейное $GL(V_*)$-представление в левом $A$-векторном пространстве $V_*$ на множество базисов левого $A$-векторного пространства $V$. Мы будем называть преобразование этого правостороннего представления на множестве базисов левого $A$-векторного пространства $V$ **активным преобразованием** потому, что гомоморфизм левого $A$-векторного пространства породил это преобразование (Смотри также определение в разделе [2]-14.1-3, а также определение на странице [1]-214).

Соответственно определению мы будем записывать действие активного преобразования $a \in GL(V_*)$ на базис $\overline{\overline{e}}$ в форме $\overline{\overline{e}}{}^*_{\,*}a$. Рассмотрим равенство

(3.32) $\qquad v^*{}_*\overline{\overline{e}}{}^*_{\,*}a = v^*{}_*\overline{\overline{e}}{}^*_{\,*}a$

Соответственно определению мы будем записывать действие активного преобразования $a \in GL(V_*)$ на базис $\overline{\overline{e}}$ в форме $\overline{\overline{e}}_*{}^*a$. Рассмотрим равенство

(3.33) $\qquad v_*{}^*\overline{\overline{e}}_*{}^*a = v_*{}^*\overline{\overline{e}}_*{}^*a$

Выражение $\overline{\overline{e}}^*{}_*a$ в левой части равенства (3.32) является образом базиса $\overline{\overline{e}}$ при активном преобразовании $a$. Выражение $v^*{}_*\overline{\overline{e}}$ в правой части равенства (3.32) является разложением вектора $\overline{v}$ относительно базиса $\overline{\overline{e}}$. Следовательно, выражение в правой части равенства (3.32) является образом вектора $\overline{v}$ относительно эндоморфизма $a$ и выражение в левой части равенства

(3.32) $\quad v^*{}_*\overline{\overline{e}}^*{}_*a = v^*{}_*\overline{\overline{e}}^*{}_*a$

является разложением образа вектора $\overline{v}$ относительно образа базиса $\overline{\overline{e}}$. Следовательно, из равенства (3.32) следует, что эндоморфизм $a$ левого $A$-векторного пространства и соответствующее активное преобразование $a$ действуют синхронно и координаты $a \circ \overline{v}$ образа вектора $\overline{v}$ относительно образа $\overline{\overline{e}}^*{}_*a$ базиса $\overline{\overline{e}}$ совпадают с координатами вектора $\overline{v}$ относительно базиса $\overline{\overline{e}}$.

Выражение $\overline{\overline{e}}_*{}^*a$ в левой части равенства (3.33) является образом базиса $\overline{\overline{e}}$ при активном преобразовании $a$. Выражение $v_*{}^*\overline{\overline{e}}$ в правой части равенства (3.33) является разложением вектора $\overline{v}$ относительно базиса $\overline{\overline{e}}$. Следовательно, выражение в правой части равенства (3.33) является образом вектора $\overline{v}$ относительно эндоморфизма $a$ и выражение в левой части равенства

(3.33) $\quad v_*{}^*\overline{\overline{e}}_*{}^*a = v_*{}^*\overline{\overline{e}}_*{}^*a$

является разложением образа вектора $\overline{v}$ относительно образа базиса $\overline{\overline{e}}$. Следовательно, из равенства (3.33) следует, что эндоморфизм $a$ левого $A$-векторного пространства и соответствующее активное преобразование $a$ действуют синхронно и координаты $a \circ \overline{v}$ образа вектора $\overline{v}$ относительно образа $\overline{\overline{e}}_*{}^*a$ базиса $\overline{\overline{e}}$ совпадают с координатами вектора $\overline{v}$ относительно базиса $\overline{\overline{e}}$.

## 4. Пассивное преобразование

Если в левом $A$-векторном пространстве $V$ определена дополнительная структура, не всякий автоморфизм сохраняет свойства заданной структуры. Например, если мы определим норму в левом $A$-векторном пространстве $V$, то для нас интересны автоморфизмы, которые сохраняют норму вектора.

**Определение 4.1.** *Нормальная подгруппа $G(V)$ группы $GL(V)$, которая порождает автоморфизмы, сохраняющие свойства заданной структуры называется* **группой симметрии**.

*Не нарушая общности, мы будем отождествлять элемент $g$ группы $G(V)$ с соответствующим автоморфизмом и записывать его действие на вектор $v \in V$ в виде $v^*{}_*g$.* □

**Определение 4.2.** *Нормальная подгруппа $G(V)$ группы $GL(V)$, которая порождает автоморфизмы, сохраняющие свойства заданной структуры называется* **группой симметрии**.

*Не нарушая общности, мы будем отождествлять элемент $g$ группы $G(V)$ с соответствующим автоморфизмом и записывать его действие на вектор $v \in V$ в виде $v_*{}^*g$.* □

**Определение 4.3.** *Мы будем называть правостороннее представление группы $G(V)$ на множестве базисов левого $A$-векторного пространства $V$* **активным левосторонним $G$-представлением**. □

Если базис $\overline{\overline{e}}$ задан, то мы можем отождествить автоморфизм левого $A$-векторного пространства $V$ и его координаты относительно базиса $\overline{\overline{e}}$. Множество $G(V_*)$ координат автоморфизмов относительно базиса $\overline{\overline{e}}$ является группой, изоморфной группе $G(V)$.

Не всякие два базиса могут быть связаны преобразованием группы симметрии потому, что не всякое невырожденное линейное преобразование принадлежит представлению группы $G(V)$. Таким образом, множество базисов можно представить как объединение орбит группы $G(V)$. В частности, если базис $\overline{\overline{e}} \in G(V)$, то орбита базиса $\overline{\overline{e}}$ совпадает с группой $G(V)$.

**Определение 4.4.** *Мы будем называть орбиту $\overline{\overline{e}}*_*G(V)$ выбранного базиса $\overline{\overline{e}}$ **многообразием базисов** левого $A$-векторного пространства $V$ столбцов.* □

**Определение 4.5.** *Мы будем называть орбиту $\overline{\overline{e}}_**G(V)$ выбранного базиса $\overline{\overline{e}}$ **многообразием базисов** левого $A$-векторного пространства $V$ строк.* □

**Теорема 4.6.** *Активное правостороннее $G(V)$-представление на многообразии базисов однотранзитивно.*

Доказательство. Теорема является следствием теоремы [3]-12.2.12. □

**Теорема 4.7.** *На многообразии базисов $\overline{\overline{e}}*_*G(V)$ левого $A$-векторного пространства столбцов существует однотранзитивное левостороннее $G(V)$-представление, перестановочное с активным.*

**Теорема 4.8.** *На многообразии базисов $\overline{\overline{e}}_**G(V)$ левого $A$-векторного пространства строк существует однотранзитивное левостороннее $G(V)$-представление, перестановочное с активным.*

Доказательство. Теорема является следствием теоремы [3]-12.2.13. □

Доказательство. Теорема является следствием теоремы [3]-12.2.14. □

Преобразование левостороннего $G(V)$-представления отличается от активного преобразования и не может быть сведено к преобразованию пространства $V$.

**Определение 4.9.** *Преобразование левостороннего $G(V)$-представления называется **пассивным преобразованием** многообразия базисов $\overline{\overline{e}}*_*G(V)$ левого $A$-векторного пространства столбцов, а левостороннее $G(V)$-представление называется **пассивным левосторонним $G(V)$-представлением**. Согласно определению мы будем записывать пассивное преобразование базиса $\overline{\overline{e}}$, порождённое элементом $a \in G(V)$, в форме $a*_*\overline{\overline{e}}$.* □

**Определение 4.10.** *Преобразование левостороннего $G(V)$-представления называется **пассивным преобразованием** многообразия базисов $\overline{\overline{e}}_**G(V)$ левого $A$-векторного пространства строк, а левостороннее $G(V)$-представление называется **пассивным левосторонним $G(V)$-представлением**. Согласно определению мы будем записывать пассивное преобразование базиса $\overline{\overline{e}}$, порождённое элементом $a \in G(V)$, в форме $a_**\overline{\overline{e}}$.* □

## 5. Геометрический объект

Активное преобразование изменяет базисы и векторы согласовано и координаты вектора относительно базиса не меняются. Пассивное преобразование меняет только базис, и это ведёт к преобразованию координат вектора относительно базиса.

Мы рассмотрим преобразование координат вектора в теореме 5.1.

Мы рассмотрим преобразование координат вектора в теореме 5.2.

**Теорема 5.1.** *Пусть $V$ - левое $A$-векторное пространство столбцов. Пусть пассивное преобразование $g \in G(V_*)$ отображает базис $\overline{\overline{e}}_1$ в базис $\overline{\overline{e}}_2$*

$$(5.1) \qquad \overline{\overline{e}}_2 = g^*{}_* \overline{\overline{e}}_1$$

*Пусть*

$$(5.2) \qquad v_i = \begin{pmatrix} v_i^1 \\ ... \\ v_i^n \end{pmatrix}$$

*матрица координат вектора $\overline{v}$ относительно базиса $\overline{\overline{e}}_i$, $i = 1, 2$. Преобразования координат*

$$(5.3) \qquad v_1 = v_2{}^*{}_* g$$

$$(5.4) \qquad v_2 = v_1{}^*{}_* g^{-1^*{}^*}$$

*не зависят от вектора $\overline{v}$ или базиса $\overline{\overline{e}}$, а определенно исключительно координатами вектора $\overline{v}$ относительно базиса $\overline{\overline{e}}$.*

**Теорема 5.2.** *Пусть $V$ - левое $A$-векторное пространство строк. Пусть пассивное преобразование $g \in G(V_*)$ отображает базис $\overline{\overline{e}}_1$ в базис $\overline{\overline{e}}_2$*

$$(5.5) \qquad \overline{\overline{e}}_2 = g_*{}^* \overline{\overline{e}}_1$$

*Пусть*

$$v_i = \begin{pmatrix} v_{i1} & ... & v_{in} \end{pmatrix}$$

*матрица координат вектора $\overline{v}$ относительно базиса $\overline{\overline{e}}_i$, $i = 1, 2$. Преобразования координат*

$$(5.6) \qquad v_1 = v_2{}_*{}^* g$$

$$(5.7) \qquad v_2 = v_1{}_*{}^* g^{-1_*{}^*}$$

*не зависят от вектора $\overline{v}$ или базиса $\overline{\overline{e}}$, а определенно исключительно координатами вектора $\overline{v}$ относительно базиса $\overline{\overline{e}}$.*

Доказательство. Теорема является следствием теоремы [3]-12.3.1. □

Предположим, что $V$, $W$ - левые $A$-векторные пространства столбцов. Пусть $G(V_*)$ - группа симметрий левого $A$-векторного пространства $V$. Гомоморфизм

$$(5.8) \qquad F : G(V_*) \to GL(W_*)$$

отображает пассивное преобразование $g \in G(V_*)$

$$(5.9) \qquad \overline{\overline{e}}_{V2} = g^*{}_* \overline{\overline{e}}_{V1}$$

левого $A$-векторного пространства $V$ в пассивное преобразование $F(g) \in$

Доказательство. Теорема является следствием теоремы [3]-12.3.2. □

Предположим, что $V$, $W$ - левые $A$-векторные пространства строк. Пусть $G(V_*)$ - группа симметрий левого $A$-векторного пространства $V$. Гомоморфизм

$$(5.13) \qquad F : G(V_*) \to GL(W_*)$$

отображает пассивное преобразование $g \in G(V_*)$

$$(5.14) \qquad \overline{\overline{e}}_{V2} = g_*{}^* \overline{\overline{e}}_{V1}$$

левого $A$-векторного пространства $V$ в пассивное преобразование $F(g) \in$

$GL(W_*)$

(5.10)  $\overline{\overline{e}}_{W2} = F(g)^*{}_*\overline{\overline{e}}_{W1}$

левого $A$-векторного пространства $W$.

Тогда координатное преобразование в левом $A$-векторном пространстве $W$ принимает вид

(5.11)  $w_2 = w_1{}^*{}_*F(g)^{-1^*{}_*}$

Следовательно, отображение

(5.12)  $F_1(g) = F(g)^{-1^*{}_*}$

является правосторонним представлением

$$F_1 : G(V_*) \longrightarrow\!\!\!*\!\!\longrightarrow W_*$$

группы $G(V_*)$ на множестве $W_*$.

**Определение 5.3.** *Мы будем называть орбиту*

(5.18)  $\mathcal{O}(V, W, \overline{\overline{e}}_V, w)$
$= (w^*{}_*F(G)^{-1^*{}_*}, G^*{}_*\overline{\overline{e}}_V)$

*представления $F_1$ многообразием* **координат геометрического объекта** *в левом $A$-векторном пространстве $V$ столбцов. Для любого базиса*

$(5.9)\quad \overline{\overline{e}}_{V2} = g^*{}_*\overline{\overline{e}}_{V1}$

*соответствующая точка*

$(5.11)\quad w_2 = w_1{}^*{}_*F(g)^{-1^*{}_*}$

*орбиты определяет* **координаты геометрического объекта** *в координатном левом $A$-векторном пространстве относительно базиса $\overline{\overline{e}}_{V2}$.*  □

**Определение 5.5.** *Допустим даны координаты $w_1$ вектора $\overline{w}$ относительно базиса $\overline{\overline{e}}_{W1}$. Мы будем называть множество векторов*

(5.20)  $\overline{w}_2 = w_2{}^*{}_*e_{W2}$

**геометрическим объектом**, *определённым в левом $A$-векторном пространстве $V$ столбцов. Для любого базиса $\overline{\overline{e}}_{W2}$, соответствующая точка*

$(5.11)\quad w_2 = w_1{}^*{}_*F(g)^{-1^*{}_*}$

$GL(W_*)$

(5.15)  $\overline{\overline{e}}_{W2} = F(g)_*{}^*\overline{\overline{e}}_{W1}$

левого $A$-векторного пространства $W$.

Тогда координатное преобразование в левом $A$-векторном пространстве $W$ принимает вид

(5.16)  $w_2 = w_{1*}{}^*F(g)^{-1}{}_*{}^*$

Следовательно, отображение

(5.17)  $F_1(g) = F(g)^{-1}{}_*{}^*$

является правосторонним представлением

$$F_1 : G(V_*) \longrightarrow\!\!\!*\!\!\longrightarrow W_*$$

группы $G(V_*)$ на множестве $W_*$.

**Определение 5.4.** *Мы будем называть орбиту*

(5.19)  $\mathcal{O}(V, W, \overline{\overline{e}}_V, w)$
$= (w_*{}^*F(G)^{-1}{}_*{}^*, G_*{}^*\overline{\overline{e}}_V)$

*представления $F_1$ многообразием* **координат геометрического объекта** *в левом $A$-векторном пространстве $V$ строк. Для любого базиса*

$(5.14)\quad \overline{\overline{e}}_{V2} = g_*{}^*\overline{\overline{e}}_{V1}$

*соответствующая точка*

$(5.16)\quad w_2 = w_{1*}{}^*F(g)^{-1}{}_*{}^*$

*орбиты определяет* **координаты геометрического объекта** *в координатном левом $A$-векторном пространстве относительно базиса $\overline{\overline{e}}_{V2}$.*  □

**Определение 5.6.** *Допустим даны координаты $w_1$ вектора $\overline{w}$ относительно базиса $\overline{\overline{e}}_{W1}$. Мы будем называть множество векторов*

(5.22)  $\overline{w}_2 = w_{2*}{}^*e_{W2}$

**геометрическим объектом**, *определённым в левом $A$-векторном пространстве $V$ строк. Для любого базиса $\overline{\overline{e}}_{W2}$, соответствующая точка*

$(5.16)\quad w_2 = w_{1*}{}^*F(g)^{-1}{}_*{}^*$

*многообразия координат определяет*

многообразия координат определяет вектор

(5.21)   $\overline{w}_2 = w_2{}^*{}_* e_{W2}$

который называется **представителем геометрического объекта** в левом $A$-векторном пространстве $V$ в базисе $\overline{\overline{e}}_{V2}$.   □

вектор

(5.23)   $\overline{w}_2 = w_{2*}{}^* e_{W2}$

который называется **представителем геометрического объекта** в левом $A$-векторном пространстве $V$ в базисе $\overline{\overline{e}}_{V2}$.   □

**Теорема 5.7. (Принцип ковариантности).** *Представитель геометрического объекта не зависит от выбора базиса $\overline{\overline{e}}_{V2}$.*

## 6. Примеры геометрического объекта

Согласно теореме 5.1, вектор является геометрическим объектом. Множество эндоморфизмов левого $A$-векторного пространства, также как и множество линейных или полилинейных отображений являются $A$-векторными пространствами. Поэтому мы можем спросить является ли эндоморфизм или полилинейное отображение геометрическим объектом.

Мы рассмотрим преобразование координат эндоморфизма левого $A$-векторного пространства столбцов в теореме 6.1.

Мы рассмотрим преобразование координат эндоморфизма левого $A$-векторного пространства строк в теореме 6.2.

**Теорема 6.1.** *Пусть $V$ - левое $A$-векторное пространство столбцов и $\overline{\overline{e}}_1$, $\overline{\overline{e}}_2$ - базисы в левом $A$-векторном пространстве $V$. Пусть базисы $\overline{\overline{e}}_1$, $\overline{\overline{e}}_2$ связаны пассивным преобразованием $g$*
$$\overline{\overline{e}}_2 = g^*{}_* \overline{\overline{e}}_1$$
*Пусть $f$ - эндоморфизм левого $A$-векторного пространства $V$. Пусть $f_i$, $i = 1, 2$, - матрица эндоморфизма $f$ относительно базиса $\overline{\overline{e}}_i$. Тогда*

(6.1)   $f_2 = g^*{}_* f_1{}^*{}_* g^{-1^*{}_*}$

**Теорема 6.2.** *Пусть $V$ - левое $A$-векторное пространство строк и $\overline{\overline{e}}_1$, $\overline{\overline{e}}_2$ - базисы в левом $A$-векторном пространстве $V$. Пусть базисы $\overline{\overline{e}}_1$, $\overline{\overline{e}}_2$ связаны пассивным преобразованием $g$*
$$\overline{\overline{e}}_2 = g_*{}^* \overline{\overline{e}}_1$$
*Пусть $f$ - эндоморфизм левого $A$-векторного пространства $V$. Пусть $f_i$, $i = 1, 2$, - матрица эндоморфизма $f$ относительно базиса $\overline{\overline{e}}_i$. Тогда*

(6.2)   $f_2 = g_*{}^* f_{1*}{}^* g^{-1}{}_*{}^*$

Доказательство. Теорема является следствием теоремы [3]-12.3.5.   □

Доказательство. Теорема является следствием теоремы [3]-12.3.6.   □

Согласно теореме 6.1, эндоморфизм векторного пространства является геометрическим объектом. Действительно, при замене базиса мы видим преобразование координат вектора и координат эндоморфизма. Но изменение координат согласовано, и образ вектора при заданом эндоморфизме не зависит от выбора базиса.

Пусть $\overline{\overline{e}}$ - базис левого векторного пространства $V$ столбцов. Линейное отображение
$$a : V \to V$$

имеет вид

$$(6.3) \quad \begin{pmatrix} v^1 \\ ... \\ v^n \end{pmatrix} = \begin{pmatrix} a_1^1 & ... & a_n^1 \\ ... & ... & ... \\ a_1^n & ... & a_n^n \end{pmatrix} \circ^\circ \begin{pmatrix} w^1 \\ ... \\ w^n \end{pmatrix} = \begin{pmatrix} a_i^1 \circ w^i \\ ... \\ a_i^n \circ w^i \end{pmatrix}$$

$$a = \begin{pmatrix} a_1^1 & ... & a_n^1 \\ ... & ... & ... \\ a_1^n & ... & a_n^n \end{pmatrix}$$

относительно базиса $\overline{\overline{e}}$. Так как частное отображение $a_j^i$ является линейным отображением $D$-алгебры $A$, то мы можем записать линейное отображение $a_j^i$ в виде

$$(6.4) \quad a_j^i = a_{j\,s0}^i \otimes a_{j\,s1}^i$$

Равенство

$$(6.5) \quad v^i = (a_{j\,s1}^i \otimes a_{j\,s2}^i) \circ w^j = a_{j\,s1}^i w^j a_{j\,s2}^i$$

является следствием равенств

$$(6.3) \quad \begin{pmatrix} v^1 \\ ... \\ v^n \end{pmatrix} = \begin{pmatrix} a_1^1 & ... & a_n^1 \\ ... & ... & ... \\ a_1^n & ... & a_n^n \end{pmatrix} \circ^\circ \begin{pmatrix} w^1 \\ ... \\ w^n \end{pmatrix} = \begin{pmatrix} a_i^1 \circ w^i \\ ... \\ a_i^n \circ w^i \end{pmatrix}$$

$$(6.4) \quad a_j^i = a_{j\,s0}^i \otimes a_{j\,s1}^i$$

**Сводка результатов 6.3.** *Пусть пассивное преобразование $g$ отображает базис $\overline{\overline{e}}_1$ в базис $\overline{\overline{e}}_2$:* $(6.6) \quad \overline{\overline{e}}_2 = g^*{}_* \overline{\overline{e}}_1$

*Пусть* $(6.7) \quad a_{k\,j}^{\ i} = a_{k\,j\,s0}^{\ i} \otimes a_{k\,j\,s1}^{\ i}$ *- координаты линейного отображения $a$ относительно базиса $\overline{\overline{e}}_k$, $k = 1, 2$. Тогда*

$$(6.13) \quad a_{2\,k\,t0}^{\ l} \otimes a_{2\,k\,t1}^{\ l} = (a_{1\,j\,s0}^{\ i} \otimes g_k^j a_{1\,j\,s1}^{\ i}) g^{-1*}{}_*{}_i^l$$

Равенства

$(6.1) \quad f_2 = g^*{}_* f_1{}^*{}_* g^{-1*}{}_*$

и

$(6.13) \quad a_{2\,k\,t0}^{\ l} \otimes a_{2\,k\,t1}^{\ l} = (a_{1\,j\,s0}^{\ i} \otimes g_k^j a_{1\,j\,s1}^{\ i}) g^{-1*}{}_*{}_i^l$

похожи. Из равенства

$$(6.14) \quad \begin{aligned} a \circ v &= \underline{a_{2\,k\,t0}^{\ l} v_2^k a_{2\,k\,t1}^{\ l} e_{2l}} = (a_{1\,j\,s0}^{\ i} v_2^k g_k^j a_{1\,j\,s1}^{\ i}) g^{-1*}{}_*{}_i^l g_l^p e_{1p} \\ &= \underline{(a_{1\,j\,s0}^{\ i} v_1^j a_{1\,j\,s1}^{\ i}) e_{1i}} \end{aligned}$$

*следует, что образ линейного отображения не зависит от выбора базиса.*

*Пусть множество тензоров* $(6.19) \quad a_{k\,j_1...j_n\,s0}^{\quad\ i} \otimes ... \otimes a_{k\,j_1...j_n\,sn}^{\quad\ i}$ *является координатами полилинейного отображения $a$ относительно базиса $\overline{\overline{e}}_k$, $k = 1, 2$. Тогда*

$$(6.23) \quad \begin{aligned} & a_2{}^l_{k_1...k_n t 0} \otimes a_2{}^l_{k_1...k_n t 1}... \otimes a_2{}^l_{k_1...k_n t n} \\ & = (a_1{}^i_{j_1...j_n s 0} \otimes g^{j_1}_{k_1} a_1{}^i_{j_1...j_n s 1}... \otimes g^{j_n}_{k_n} a_1{}^i_{j_1...j_n s n}) g^{-1*}{}^l_i \end{aligned}$$

*Из равенства*

$$(6.24) \quad \begin{aligned} a \circ (v_1,...,v_2) & = a_2{}^l_{k_1...k_n t 0} v^{k_1}_{12} a_2{}^l_{k_1...k_n t 1}...v^{k_n}_{n 2} a_2{}^l_{k_1...k_n t n} e_{2l} \\ & = (a_1{}^i_{j_1...j_n s 0} v^{k_1}_{12} g^{j_1}_{k_1} a_1{}^i_{j_1...j_n s 1}...v^{k_n}_{n 2} g^{j_n}_{k_n} a_1{}^i_{j_1...j_n s n}) g^{-1*}{}^l_i g^p_l e_{1p} \\ & = (a_1{}^i_{j_1...j_n s 0} v^{j_1}_{11} a_1{}^i_{j_1...j_n s 1}...v^{j_n}_{n 1} a_1{}^i_{j_1...j_n s n}) e_{1i} \end{aligned}$$

*следует, что образ полилинейного отображения не зависит от выбора базиса.*

*Рассмотрим кососимметричное полилинейное отображение*

$$(6.26) \quad \begin{aligned} h^i \circ (u \otimes v) & = a^i_{jk s 0} u^j a^i_{jk s 1} v^k a^i_{jk s 2} \\ & = (a^i_{jk s 0} \otimes a^i_{jk s 1} \otimes a^i_{jk s 2}) \circ (u^j \otimes v^k) \end{aligned}$$

*Если мы рассмотрим выражение, подобное определителю,*

$$(6.28) \quad \det{}^* \begin{pmatrix} u^j & v^j \\ u^k & v^k \end{pmatrix} = u^j \otimes v^k - v^j \otimes u^k$$

*то мы можем записать полилинейное отображение $h$ в виде*

$$(6.29) \quad h^i \circ (u \otimes v) = \tfrac{1}{2} (a^i_{jk s 0} \otimes a^i_{jk s 1} \otimes a^i_{jk s 2}) \circ \det{}^* \begin{pmatrix} u^j & v^j \\ u^k & v^k \end{pmatrix}$$

*Пусть множество тензоров* $(6.19) \quad a^i_{k j_1...j_n s 0} \otimes ... \otimes a^i_{k j_1...j_n s n}$ *является координатами полилинейного отображения $h$ относительно базиса $\bar{\bar{e}}_k$, $k = 1, 2$. Тогда*

$$(6.37) \quad \begin{aligned} & (a_2{}^i_{jk s 0} \otimes_1 a_2{}^i_{jk s 1} \otimes_2 a_2{}^i_{jk s 2}) \circ \det{}^* \begin{pmatrix} 1 \otimes_1 \delta^j_q & 1 \otimes_2 \delta^j_r \\ 1 \otimes_1 \delta^k_q & 1 \otimes_2 \delta^k_r \end{pmatrix} \\ & = ((a_1{}^i_{jk s 0} \otimes_1 a_1{}^i_{jk s 1} \otimes_2 a_1{}^i_{jk s 2}) g^{-1*}{}^p_i) \circ \det{}^* \begin{pmatrix} 1 \otimes_1 g^j_q & 1 \otimes_2 g^j_r \\ 1 \otimes_1 g^k_q & 1 \otimes_2 g^k_r \end{pmatrix} \end{aligned}$$

$\square$

Пусть пассивное преобразование $g$ отображает базис $\bar{\bar{e}}_1$ в базис $\bar{\bar{e}}_2$

$$(6.6) \quad \bar{\bar{e}}_2 = g^*{}_* \bar{\bar{e}}_1$$

Пусть

$$(6.7) \quad a^i_{k j} = a^i_{k j s 0} \otimes a^i_{k j s 1}$$

координаты линейного отображения $a$ относительно базиса $\bar{\bar{e}}_k$, $k = 1, 2$.

Пусть $v^j_k$, $w^j_k$ координаты векторов $v$, $w$ относительно базиса $\bar{\bar{e}}_k$, $k = 1, 2$. Тогда

$$(6.8) \quad w^i_k = a^i_{k j s 0} v^j_k a^i_{k j s 1}$$

Согласно теореме 6.1,

$$v_1 = v_2{}^*{}_*g \tag{6.9}$$

$$w_1 = w_2{}^*{}_*g \tag{6.10}$$

Равенства

$$w_2^l g_l^i = a_1{}^i_{j\,s0} v_2^k g_k^j a_1{}^i_{j\,s1} \tag{6.11}$$

$$w_2^l = (a_1{}^i_{j\,s0} v_2^k g_k^j a_1{}^i_{j\,s1}) g^{-1}{}^*{}_*{}^l_i \tag{6.12}$$

являются следствием равенств

$\boxed{(6.9) \quad v_1 = v_2{}^*{}_*g}$ $\boxed{(6.10) \quad w_1 = w_2{}^*{}_*g}$

$\boxed{(6.8) \quad w_1^i = a_1{}^i_{j\,s0} v_1^j a_1{}^i_{j\,s1}}$ (k=1).

Равенство

$$a_2{}^l_{k\,t0} \otimes a_2{}^l_{k\,t1} = (a_1{}^i_{j\,s0} \otimes g_k^j a_1{}^i_{j\,s1}) g^{-1}{}^*{}_*{}^l_i \tag{6.13}$$

является следствием равенств

$\boxed{(6.8) \quad w_2^i = a_2{}^i_{j\,s0} v_2^j a_2{}^i_{j\,s1}}$ (k=2),

$\boxed{(6.12) \quad w_2^l = (a_1{}^i_{j\,s0} v_2^k g_k^j a_1{}^i_{j\,s1}) g^{-1}{}^*{}_*{}^l_i}$

Равенство

$$a \circ v = \underline{a_2{}^l_{k\,t0} v_2^k a_2{}^l_{k\,t1} e_{2l}} = (a_1{}^i_{j\,s0} v_2^k g_k^j a_1{}^i_{j\,s1}) g^{-1}{}^*{}_*{}^l_i g_l^p e_{1p} \tag{6.14}$$
$$= \underline{(a_1{}^i_{j\,s0} v_1^j a_1{}^i_{j\,s1}) e_{1i}}$$

является следствием равенств

$\boxed{(6.9) \quad v_1 = v_2{}^*{}_*g}$

$\boxed{(6.13) \quad a_2{}^l_{k\,t0} \otimes a_2{}^l_{k\,t1} = (a_1{}^i_{j\,s0} \otimes g_k^j a_1{}^i_{j\,s1}) g^{-1}{}^*{}_*{}^l_i}$

Из равенства (6.14) следует, что образ линейного отображения не зависит от выбора базиса.

Нетрудно видеть сходство равенств

$\boxed{(6.1) \quad f_2 = g^*{}_* f_1{}^*{}_* g^{-1}{}^*{}_*}$

и

$\boxed{(6.13) \quad a_2{}^l_{k\,t0} \otimes a_2{}^l_{k\,t1} = (a_1{}^i_{j\,s0} \otimes g_k^j a_1{}^i_{j\,s1}) g^{-1}{}^*{}_*{}^l_i}$

Если $D$-алгебра $A$ коммутативна, эти равенства совпадают.

**Теорема 6.4.** *Полилинейное отображение*

$$a : V^n \to V \tag{6.15}$$

*является геометрическим объектом.*

Доказательство. Не нарушая общности, мы рассмотрим доказательство теоремы для левого $A$-векторного пространства столбцов.

Пусть пассивное преобразование $g$ отображает базис $\overline{\overline{e}}_1$ в базис $\overline{\overline{e}}_2$

$$\overline{\overline{e}}_2 = g^*{}_* \overline{\overline{e}}_1 \tag{6.16}$$

Пусть $v^j_{ik}$, $w^j_k$ координаты векторов $v_i$, $i = 1, ..., n$, $w$ относительно базиса $\overline{\overline{e}}_k$, $k = 1, 2$. Согласно теореме 6.1,

$$(6.17) \qquad v_{i1} = v_{i2} {}^*{}_* g$$

$$(6.18) \qquad w_1 = w_2 {}^*{}_* g$$

Пусть множество тензоров

$$(6.19) \qquad a^{\,i}_{k\,j_1...j_n}{}_{s0} \otimes ... \otimes a^{\,i}_{k\,j_1...j_n}{}_{sn}$$

является координатами полилинейного отображения $a$ относительно базиса $\overline{\overline{e}}_k$. Тогда[6.1]

$$(6.20) \qquad w^i_k = a^{\,i}_{k\,j_1...j_n}{}_{s0} v^{j_1}_{1k} a^{\,i}_{k\,j_1...j_n}{}_{s1} ... v^{j_n}_{nk} a^{\,i}_{k\,j_1...j_n}{}_{sn}$$

Равенства

$$(6.21) \qquad w^l_2 g^i_l = a^{\,i}_{1\,j_1...j_n}{}_{s0} v^{k_1}_{12} g^{j_1}_{k_1} a^{\,i}_{1\,j_1...j_n}{}_{s1} ... v^{k_n}_{n2} g^{j_n}_{k_n} a^{\,i}_{k\,j_1...j_n}{}_{sn}$$

$$(6.22) \qquad w^l_2 = (a^{\,i}_{1\,j_1...j_n}{}_{s0} v^{k_1}_{12} g^{j_1}_{k_1} a^{\,i}_{1\,j_1...j_n}{}_{s1} ... v^{k_n}_{n2} g^{j_n}_{k_n} a^{\,i}_{k\,j_1...j_n}{}_{sn}) g^{-1\,*\,l}{}_*{}_i$$

являются следствием равенств

$$\boxed{(6.17) \quad v_{i1} = v_{i2} {}^*{}_* g} \qquad \boxed{(6.18) \quad w_1 = w_2 {}^*{}_* g}$$

$$\boxed{(6.20) \quad w^i_1 = a^{\,i}_{1\,j_1...j_n}{}_{s0} v^{j_1}_{11} a^{\,i}_{1\,j_1...j_n}{}_{s1} ... v^{j_n}_{n1} a^{\,i}_{1\,j_1...j_n}{}_{sn}} \qquad (k=1).$$

Равенство

$$(6.23) \quad \begin{aligned} & a^{\,l}_{2\,k_1...k_n}{}_{t0} \otimes a^{\,l}_{2\,k_1...k_n}{}_{t1} ... \otimes a^{\,l}_{2\,k_1...k_n}{}_{tn} \\ &= (a^{\,i}_{1\,j_1...j_n}{}_{s0} \otimes g^{j_1}_{k_1} a^{\,i}_{1\,j_1...j_n}{}_{s1} ... \otimes g^{j_n}_{k_n} a^{\,i}_{1\,j_1...j_n}{}_{sn}) g^{-1\,*\,l}{}_*{}_i \end{aligned}$$

является следствием равенств

$$\boxed{(6.20) \quad w^i_2 = a^{\,i}_{2\,j_1...j_n}{}_{s0} v^{j_1}_{12} a^{\,i}_{2\,j_1...j_n}{}_{s1} ... v^{j_n}_{n2} a^{\,i}_{2\,j_1...j_n}{}_{sn}} \qquad (k=2),$$

$$\boxed{(6.22) \quad w^l_2 = (a^{\,i}_{1\,j_1...j_n}{}_{s0} v^{k_1}_{12} g^{j_1}_{k_1} a^{\,i}_{1\,j_1...j_n}{}_{s1} ... v^{k_n}_{n2} g^{j_n}_{k_n} a^{\,i}_{k\,j_1...j_n}{}_{sn}) g^{-1\,*\,l}{}_*{}_i}$$

Равенство

$$(6.24) \quad \begin{aligned} a \circ (v_1, ..., v_2) &= \underline{a^{\,l}_{2\,k_1...k_n}{}_{t0} v^{k_1}_{12} a^{\,l}_{2\,k_1...k_n}{}_{t1} ... v^{k_n}_{n2} a^{\,l}_{2\,k_1...k_n}{}_{tn} e_{2l}} \\ &= (a^{\,i}_{1\,j_1...j_n}{}_{s0} v^{k_1}_{12} g^{j_1}_{k_1} a^{\,i}_{1\,j_1...j_n}{}_{s1} ... v^{k_n}_{n2} g^{j_n}_{k_n} a^{\,i}_{1\,j_1...j_n}{}_{sn}) g^{-1\,*\,l}{}_*{}_i g^p_l e_{1p} \\ &= \underline{(a^{\,i}_{1\,j_1...j_n}{}_{s0} v^{j_1}_{11} a^{\,i}_{1\,j_1...j_n}{}_{s1} ... v^{j_n}_{n1} a^{\,i}_{1\,j_1...j_n}{}_{sn}) e_{1i}} \end{aligned}$$

является следствием равенств

$$\boxed{(6.17) \quad v_{i1} = v_{i2} {}^*{}_* g}$$

$$\boxed{(6.23) \quad \begin{aligned} & a^{\,l}_{2\,k_1...k_n}{}_{t0} \otimes a^{\,l}_{2\,k_1...k_n}{}_{t1} ... \otimes a^{\,l}_{2\,k_1...k_n}{}_{tn} \\ &= (a^{\,i}_{1\,j_1...j_n}{}_{s0} \otimes g^{j_1}_{k_1} a^{\,i}_{1\,j_1...j_n}{}_{s1} ... \otimes g^{j_n}_{k_n} a^{\,i}_{1\,j_1...j_n}{}_{sn}) g^{-1\,*\,l}{}_*{}_i \end{aligned}}$$

Из равенства (6.24) следует, что образ полилинейного отображения не зависит от выбора базиса. $\square$

---

[6.1] В этом доказательстве нет необходимости рассматривать возможность изменений порядка переменных в каждом слагаемом, так как нас интересует одно слагаемое. Если в этом слагаемом меняется порядок переменных, то мы можем изменить этот порядок до того, как мы начнём вычисления в доказательстве.

Мы рассмотрим преобразование кососимметричного полилинейного отображения на примере билинейного отображения. Положим

$$(6.25) \qquad w = h \circ (u \otimes v) = \frac{1}{2} h \circ (u \otimes v - v \otimes u)$$

где отображение $h$ имеет вид

$$(6.26) \qquad \begin{aligned} h^i \circ (u \otimes v) &= a^i_{jk\,s0} u^j a^i_{jk\,s1} v^k a^i_{jk\,s2} \\ &= (a^i_{jk\,s0} \otimes a^i_{jk\,s1} \otimes a^i_{jk\,s2}) \circ (u^j \otimes v^k) \end{aligned}$$

Равенство

$$(6.27) \qquad w^i = h^i \circ (u \otimes v) = \frac{1}{2}(a^i_{jk\,s0} \otimes_1 a^i_{jk\,s1} \otimes_2 a^i_{jk\,s2}) \circ (u^j \otimes v^k - v^j \otimes u^k)$$

является следствием равенств (6.25), (6.26). Если мы рассмотрим выражение, подобное определителю,

$$(6.28) \qquad \det{}^* \begin{pmatrix} u^j & v^j \\ u^k & v^k \end{pmatrix} = u^j \otimes v^k - v^j \otimes u^k$$

то равенство (6.27) примет вид

$$(6.29) \qquad h^i \circ (u \otimes v) = \frac{1}{2}(a^i_{jk\,s0} \otimes a^i_{jk\,s1} \otimes a^i_{jk\,s2}) \circ \det{}^* \begin{pmatrix} u^j & v^j \\ u^k & v^k \end{pmatrix}$$

Пусть пассивное преобразование $g$ отображает базис $\overline{\overline{e}}_1$ в базис $\overline{\overline{e}}_2$

$$(6.30) \qquad \overline{\overline{e}}_2 = g^*{}_* \overline{\overline{e}}_1$$

Равенство

$$(6.31) \qquad \begin{aligned} w_2^p g^i_p &= \frac{1}{2}(a_1{}^i_{jk\,s0} \otimes_1 a_1{}^i_{jk\,s1} \otimes_2 a_1{}^i_{jk\,s2}) \\ &\quad \circ ((u_2^q g^j_q) \otimes (v_2^r g^k_r) - (v_2^r g^j_r) \otimes (u_2^q g^k_q)) \\ &= \frac{1}{2}(a_1{}^i_{jk\,s0} \otimes_1 a_1{}^i_{jk\,s1} \otimes_2 a_1{}^i_{jk\,s2}) \\ &\quad \circ ((1 \otimes_1 g^j_q) \otimes (1 \otimes_2 g^k_r) - (1 \otimes_2 g^j_r) \otimes (1 \otimes_1 g^k_q)) \circ (u_2^q \otimes v_2^r) \end{aligned}$$

является следствием равенств

$$\boxed{(6.27) \quad w^i = h^i \circ (u \otimes v) = \tfrac{1}{2}(a^i_{jk\,s0} \otimes_1 a^i_{jk\,s1} \otimes_2 a^i_{jk\,s2}) \circ (u^j \otimes v^k - v^j \otimes u^k)} \quad \boxed{(6.18) \quad w_1 = w_2{}^*{}_* g}$$

$$(6.32) \qquad u_1 = u_2{}^*{}_* g$$

$$(6.33) \qquad u_1 = u_2{}^*{}_* g$$

Равенство

$$(6.34) \qquad \begin{aligned} w_2^p &= \frac{1}{2}((a_1{}^i_{jk\,s0} \otimes_1 a_1{}^i_{jk\,s1} \otimes_2 a_1{}^i_{jk\,s2}) g^{-1}{}_*{}^{*\,p}_{\ i}) \\ &\quad \circ ((1 \otimes_1 g^j_q) \otimes (1 \otimes_2 g^k_r) - (1 \otimes_2 g^j_r) \otimes (1 \otimes_1 g^k_q)) \circ (u_2^q \otimes v_2^r) \end{aligned}$$

является следствием равенства (6.31). Равенство

$$(6.35) \qquad \begin{aligned} w_2^p &= \frac{1}{2}(a_2{}^i_{jk\,s0} \otimes_1 a_2{}^i_{jk\,s1} \otimes_2 a_2{}^i_{jk\,s2}) \\ &\quad \circ ((1 \otimes_1 \delta^j_q) \otimes (1 \otimes_2 \delta^k_r) - (1 \otimes_2 \delta^j_r) \otimes (1 \otimes_1 \delta^k_q)) \circ (u_2^q \otimes v_2^r) \end{aligned}$$

является следствием равенства

$$(6.27) \quad w^i = h^i \circ (u \otimes v) = \frac{1}{2}(a^i_{jk\,s0} \otimes_1 a^i_{jk\,s1} \otimes_2 a^i_{jk\,s2}) \circ (u^j \otimes v^k - v^j \otimes u^k)$$

Равенство

$$(6.36) \quad \begin{aligned} &(a_{2\,jk\,s0}^i \otimes_1 a_{2\,jk\,s1}^i \otimes_2 a_{2\,jk\,s2}^i) \\ &\circ \left((1 \otimes_1 \delta_q^j) \otimes (1 \otimes_2 \delta_r^k) - (1 \otimes_2 \delta_r^j) \otimes (1 \otimes_1 \delta_q^k)\right) \\ &= ((a_{1\,jk\,s0}^i \otimes_1 a_{1\,jk\,s1}^i \otimes_2 a_{1\,jk\,s2}^i) g^{-1}{_*}{^*}{_i^p}) \\ &\circ \left((1 \otimes_1 g_q^j) \otimes (1 \otimes_2 g_r^k) - (1 \otimes_2 g_r^j) \otimes (1 \otimes_1 g_q^k)\right) \end{aligned}$$

является следствием равенства (6.34), (6.35). Равенство

$$(6.37) \quad \begin{aligned} &(a_{2\,jk\,s0}^i \otimes_1 a_{2\,jk\,s1}^i \otimes_2 a_{2\,jk\,s2}^i) \circ \det{}^* \begin{pmatrix} 1 \otimes_1 \delta_q^j & 1 \otimes_2 \delta_r^j \\ 1 \otimes_1 \delta_q^k & 1 \otimes_2 \delta_r^k \end{pmatrix} \\ &= ((a_{1\,jk\,s0}^i \otimes_1 a_{1\,jk\,s1}^i \otimes_2 a_{1\,jk\,s2}^i) g^{-1}{_*}{^*}{_i^p}) \circ \det{}^* \begin{pmatrix} 1 \otimes_1 g_q^j & 1 \otimes_2 g_r^j \\ 1 \otimes_1 g_q^k & 1 \otimes_2 g_r^k \end{pmatrix} \end{aligned}$$

является следствием равенства (6.28), (6.36).

## 7. Список литературы

# 8. Предметный указатель

## 9. Специальные символы и обозначения